\documentclass[12pt]{report}
\usepackage{graphicx}
\usepackage{amsmath}
\usepackage{amscd}
\usepackage{amsthm}
\usepackage{amsfonts}
\usepackage{amssymb}
\usepackage[alphabetic,backrefs,lite]{amsrefs}
\usepackage{mathrsfs,amsrefs}
\usepackage{longtable,cases}

\newtheorem{thm}{Theorem}[section]

\newtheorem{claim}[thm]{Claim}

\theoremstyle{definition}
\newtheorem{exam}[thm]{Example}
\newtheorem{defn}[thm]{Definition}
\theoremstyle{plain}

\newtheorem{lem}[thm]{Lemma}

\newtheorem{problem}[thm]{Problem}
\newtheorem{prop}[thm]{Proposition}
\newtheorem{rem}[thm]{Remark}

\newcommand{\bF}{\mathbb F}
\newcommand{\bP}{\mathbb P}

\begin{document}
\pagenumbering{roman} \setcounter{page}{2}

\renewcommand{\title}{\hspace{.3cm} Rational Points on del Pezzo Surfaces of degree 1 and 2}
\renewcommand{\author}{Shuijing Li}


 \thispagestyle{empty}
 {\centering \section*{Abstract} }
 \addcontentsline{toc}{section}{Abstract}
 {\centering \Large
\title \par}
{\centering by \par}  {\centering \large \author
\par} \vspace{0.25in}

  One of the fundamental problems in Algebraic Geometry is to study
solutions to certain systems of polynomial equations in several variables, or in other words, find rational points on a given variety which is defined by equations. In this paper, we discuss the existence of del Pezzo surface of degree 1 and 2 with a unique rational point over any finite field $\mathbb F_q$, and we will give a lower bound on the number of rational points to each $q$. Furthermore, we will give explicit equations of del Pezzo surfaces with a unique rational point.

Also, we will discuss the rationality property of the del Pezzo surfaces especially in lower degrees.

\newpage
\tableofcontents

\newpage
 \listoffigures
\addcontentsline{toc}{section}{List of Figures}


\chapter{Introduction}\label{intro}
\pagenumbering{arabic} \setcounter{page}{1}

One of the fundamental problems in Arithmetic Geometry is to find
solutions to Diophantine equations, or equivalently, the existence of
rational points on projective algebraic varieties, over finite
fields of small cardinality. In other words, to determine or
describe the set of rational points X(\textit{k}) for a given
algebraic variety X defined over a field \textit{k}.

  \begin{defn}
Given a field $k$, a $k$-variety X is a geometrically integral, separated scheme of finite type over $k$.
  \end{defn}




\begin{exam}
Plane algebraic curves, which include lines, circles, parabolas and cubic curves form one of the best studied classes of algebraic varieties.
\end{exam}

\begin{exam}

Over the complex number $\mathbb C$, the simplest quadratic equation $x^2+ 1 = 0 $ defines a variety X. It has two roots $ x= \pm i$, that is to say, $X(\mathbb C)= \{\pm i \}$.

But $x^2+ 1 = 0 $ doesn't have any solution over $\mathbb R$, so $X(\mathbb R) = \emptyset$.

However, $x^2+ 1 = 0 $ has solution $x=\bar 1$ over the finite field $\mathbb F_2$, and has no solution again over $\mathbb F_3$. Therefore, X has no rational point over $\mathbb F_3$, and X has a unique rational point over $\mathbb F_2$.

\end{exam}

 In this chapter, we will define del Pezzo surfaces and some basic geometric properties. We will also introduce the rationality and the unirationality properties of surfaces, especially of del Pezzo surfaces.
\section{del Pezzo surfaces}
A del Pezzo surface over a field \textit{k} is a smooth
projective algebraic surface X for which the anticanonical
sheaf $\omega_X^{-1}=-K_X$ is ample. The degree of X is the integer
$d=(\omega_X^{-1},\omega_X^{-1})$ (where (, ) denotes intersection
form). Let $\bar{X}= X \times_k \bar{k}$, where $\bar{\textit{k}}$
is the algebraic closure of \textit{k}.

We recall that the property of being ample by definition means that there exists a closed embedding $\textit{i}:
X\hookrightarrow \mathbb P^N$ and some power $\emph{n}>0 $ such that $ \omega_X^{-n} \cong
\textit{i}^{*}( \textit{O}_{ \mathbb P^{N}} (1))$.

Many interesting arithmetic questions are connected with the class of del Pezzo surfaces, as
such surfaces are geometrically rational (they are rational over the
field $\bar{k}$). It is especially interesting to look at
problems concerning the existence of \textit{k}-rational points. From the arithmetic and geometric point
of view, higher degree del Pezzo surfaces are simpler than the low
degree ones (see proposition \ref{Prop}).

\begin{defn}
 A set of points in $\mathbb P^2$ is in \textbf{general position} if there are no three points on a line, no six on a conic and no eight on a singular cubic with
a singularity at one of the points.
\end{defn}

  To sum up the properties of del Pezzo surfaces that we will use, we combine
several facts taken from \cite{Manin} into the following proposition.

\begin{prop}\label{Prop}

    Let X be a del Pezzo surface of degree d over a field \textit{k}. \\
(a) The degree of X satisfies $1 \leq  d \leq  9.$ \\
(b) The Picard group Pic$\bar{X}$ is a free abelian group of rank $10 - d$. \\
(c) If $X\longrightarrow X' $ is a birational morphism then $X'$ is
a del Pezzo surface. \\
(d) Over $\bar{\textit{k}}$, either $\bar{X}$ is isomorphic to the
blowup of $\mathbb P^2 _ {\bar k} $ at $r = 9-d$ points $ \{ x_1, . . . , x_r \} $ in general
position, or $d = 8$ and $\bar{X} \cong {\mathbb P^1_{\bar{k}} }
\times  \mathbb P^1_{\bar{k}} $

 Conversely, any surface that is realized as the blow-up of a projective plane at r
points in general position is a del Pezzo surface of degree $d = 9 -
r$.

(e) The anticanonical map $X \longrightarrow \mathbb P^d$ is a
closed immersion for $d \geq  3$, so X can be realized as a degree d
surface in $\mathbb P^d$ for $d \geq  3$. The set of exceptional
curves is identified under this embedding with the set of lines in
the ambient space which lie on X. This is not true for $d=1$ or 2, where
X can be realized only in weighted projective spaces. We will give explicit equations in the following chapters.
\end{prop}
\begin{proof}
   These are (respectively) Theorem 24.3(i), Lemma 24.3.1, Remark 24.4.2, Theorem
24.4, Corollary 24.5.2, Remark 24.5.1, and Theorem 26.2 of
\cite{Manin}.
\end{proof}

\section{Exceptional curves on del Pezzo surfaces}

 Let X be a del Pezzo surface of degree d over a field $k$.
\begin{defn} For X isomorphic to the blowup of the projective plane at $r = 9- d$ points  $ \{ x_1,\cdots, x_r \} $ in
general position, a curve C is an \textbf{exceptional curve} or \textbf{$(-1)$-curve} on X if the self intersection number of C is $(-1)$.

For $d \geq 3$, the image of an exceptional curve under the anticanonical embedding map is a line (i.e. $C \cong \mathbb P^1_{\bar k}$).
\end{defn}
\begin{rem}(\cite{Manin})
 From the definition, we can get that the image of C under the blowing-down map to the projective plane is either:
\begin{enumerate}
\item one of the $x_i$
\item a line passing through two of the $x_i$
\item a conic passing through five of the $x_i$
\item a cubic passing through seven of the $x_i$ such that one $x_i$ is a double point
\item a quartic passing through eight of the $x_i$ such that three $x_i$ are double points
\item a quintic passing through eight of the $x_i$ such that six $x_i$ are double points
\item a sextic passing through eight of the $x_i$ such that seven $x_i$ are double points and
one is a triple point
\end{enumerate}
\end{rem}

We will speak of ``exceptional curves on X'' instead of $\bar{X}$. In the above proposition part (b), the Picard group has generators that can be described in terms of exceptional curves on the surface\textbf{\cite{Manin}}.

The blow down of any exceptional curve on a del Pezzo surface is a del Pezzo surface of degree 1 more. The blow up of any point on a del Pezzo surface is a del Pezzo surface of degree 1 less, provided that the point does not lie on any exceptional curve and the degree is greater than 1. Moreover, the point shouldn't lie on the branch locus when the degree is 2.

\begin{rem}
 There are only finitely many exceptional curves on del Pezzo
surfaces and these depend only on the degree (unless the degree is 8). For each value of $r= 9-d$, the number of exceptional curves is easily computable:
\begin{itemize}
 \item del Pezzo surfaces of degree 1 have 240 exceptional curves and the integral automorphism group of Pic(X) that fixes $-K_X$ is the Weyl group of the root system $E_8$. We will present a detailed discussion of these facts in chapter 3.
\item del Pezzo surfaces of degree 2 have 56 exceptional curves and the integral automorphism group of Pic(X) that fixes $-K_X$ is the Weyl group of the root system $E_7$. The linear system of its anticanonical divisor defines a map from the del Pezzo surface to the projective plane, branched over a quartic plane curve. We will give more discussion of this in chapter 4.
\item del Pezzo surfaces of degree 3 are the cubic surfaces in $\bP^3$. They have 27 exceptional curves, and the integral automorphism group of Pic(X) that fixes $-K_X$ is the Weyl group of $E_6$. We will give more discussion of this in chapter 5.
\item del Pezzo surfaces of degree 4 are given by the intersection of two quadrics in $\bP^4$. They have 16 exceptional curves.
\item del Pezzo surfaces of degree 5 have 10 exceptional curves. There is only one such surface up to isomorphism over the algebraic closure.
\item del Pezzo surfaces of degree 6 have 6 exceptional curves. There is only one such surface up to isomorphism over the algebraic closure.
\item del Pezzo surfaces of degree 7 have 3 exceptional curves. There is only one such surface up to isomorphism.
 \item del Pezzo surfaces of degree 8 have 2 isomorphism types over the algebraic closure. One is a Hirzebruch surface given by the blow up of the projective plane at one point. This surface has a unique exceptional curve. The other is ${\mathbb P^1_{\bar{k}} }\times  \mathbb P^1_{\bar{k}}$, which is the only del Pezzo surface that cannot be obtained by blowing up points on the projective plane. This one doesn't contain exceptional curves.
\end{itemize}
\end{rem}
\section{Rationality and unirationality}

Determining the rational points on a del Pezzo surface is usually very difficult, but there are some answers when we restrict to special fields $k$. Related to the rational points on X, we can look for other geometric properties of X as well. First, let's recall the definition of rational variety:
\begin{defn}
  A variety X over $k$ is \textbf{rational} if it is birational to $\bP^n_{\bar k}$, for some n.
\end{defn}

\begin{exam}
Examples of rational surfaces:
\begin{itemize}
  \item The projective plane.
  \item del Pezzo surfaces.
  \item  Veronese surface, which is an embedding of the projective plane into $\mathbb P^5$.

\end{itemize}
\end{exam}
\textbf{Note:}
When we say a variety is ``rational", we mean it is rational over $\bar k$; otherwise we need to specify the ground field $k$.
The question of whether a variety is rational is one of the most basic questions one can ask. When we work over a non-algebraically closed field $k$, this question becomes especially tricky. Various definitions are relevant, and we will focus on the property called ``unirationality".
\begin{defn}
A variety X over $k$ is \textbf{unirational} if there is a dominant rational map $\bP^n \dashrightarrow X$, defined over $k$.
 \end{defn}
\textbf{Note:}
The ground field $k$ must be specified when we consider the ``unirationality".
 Here are some basic properties of unirational varieties:
\begin{prop}
 Let X be a variety over $k$. The following are equivalent:\\
(1) X is unirational.\\
(2) There is a dominant generically finite morphism $\phi: Y \longrightarrow X$, where Y is rational over $k$.\\
(3) The function field of X is contained in a purely transcendental field extension of $k$.
\end{prop}

\begin{thm}
A smooth del Pezzo surface X of degree 2 having a $k$-point not lying on the exceptional curves, nor on the ramification curve of the map given by $|-K_X| : X \longrightarrow \mathbb P^2$, is
 $k$-unirational.
\end{thm}
\begin{proof}
 For the latter case, the proof is in \textbf{\cite{Manin}}; for degenerate ones see Coray and Tsfasman \textbf{\cite{Cm}}.
\end{proof}
\begin{rem}
Segre proved that a smooth cubic surface over $\mathbb Q$ is unirational if and only if it has a rational point. J\'anos Koll\'ar proved in 2000 that this result also holds for smooth cubic hypersurfaces of dimension at least 2 over any field, including finite fields\textbf{\cite{Kollar02}}.

\end{rem}
The question of \textit{k}-unirationality  of del Pezzo surfaces of degree 1 and
conic bundles of degree at most 1 is an extremely difficult problem.
However over local fields \textit{k}, a conic bundle having a
\textit{k}-point is always unirational (see Theorem 3.5.2 of
\cite{Yan})

There are some results showing that, over the field of rational
functions, the number of rational points on certain del Pezzo surfaces
is infinite. However, over finite fields, the number of
$k$-rational points obviously cannot be infinite anymore. Hence we turn to analyzing the minimal number of rational points. In this paper, I will give the lower bound of $k$-rational points
over all finite fields on del Pezzo surfaces of degree 1 and 2.

\chapter{General discussion of rational points on del Pezzo surfaces}

Let $F$ be a homogeneous polynomial over $\bF_q$ of degree $d$ in
$n>d$ variables; the Chevalley-Warning theorem says that $F$ has a
nontrivial zero. In geometric terms, the set of $\bF_q$-rational
points $X(\bF_q)$ of the hypersurface $X=\{F=0\}$ has positive
cardinality (in fact, $|X(\bF_q)|\equiv 1$ mod $q$).

 Let $k$ be a finite field of q elements and X be a smooth projective surface over $k$. Let G denote the Galois group $G= Gal( \bar k /k)$ and $F\in G$ denote the Frobenius endomorphism of $\bar{k}/k$ given by $$F:z\longrightarrow z^q  (z \in \bar{k}). $$
  Since the classes of the lines on $X \otimes \bar k$ are defined over $\bar{k}$ and generate the group $Pic(X \otimes \bar k)$, F induces a permutation $F^* $ on the lines. Let N denote the number of $k$-points of the surface
  X.
\begin{thm}
 \emph{(A.Weil)} If X is a smooth projective rational surface over $k= \mathbb F_q$, then
\begin{equation}
  N= |X(k)| = q^2 + q Tr F^{*}  + 1 \nonumber,
\end{equation}
where $Tr F^*$ denotes the trace of $F^*$ in the representation of
$Gal(\bar k /k)$ on $Pic(X \otimes \bar k)$.
\end{thm}
This theorem is due to Weil, as stated Theorem 23.1 in \textbf{\cite{Manin}}.

A necessary condition for the existence of a del Pezzo surface with
one rational point is the existence of a conjugacy class $\sigma=F^*$ in the Weyl group with $$ Tr \sigma  =-q.$$

\begin{exam}
A smooth cubic surface in $\mathbb P^3$ is a del Pezzo surface of
degree 3. If the cubic surface contains a unique rational point (i.e.
N=1), by Weil's Theorem,
\begin{equation}
   1= q^2 + q  Tr F^*  + 1 \nonumber
\end{equation}

By table 1 on CH IV, $\S 31$ \textbf{\cite{Manin}} (classes of conjugacy elements in the Weyl group $W(E_6)$), the smallest trace is -2, hence $q=2$ is the only possible value to satisfy the equation.

The existence and uniqueness of such surfaces with a unique rational point has been studied by many people. H.P.F. Swinnerton-Dyer \textbf{\cite{Sw4}} stated that a smooth projective cubic hypersurface defined over a finite field $k=\mathbb F_q$ with a unique rational point could only happen when $q=2$.

 Explicitly, the unique smooth cubic hypersurface with a unique rational point up to linear transformation over $\mathbb
F_2$ in $\mathbb P^3$ is given by:
$$F(x,y,z,w)=xw^{2}+x^{2}
      w+y^{3}+z^{2}y+z^{3}+x^{3}+ yzx+(y+z)x^{2} $$

\end{exam}

For other del Pezzo surfaces, we have some results as follows:

\begin{thm}
\emph (Enriques, Swinnerton-Dyer) Every del Pezzo surface X of
degree 5 defined over $k$ has a $k$-point.
\end{thm}

This theorem was first formulated by Enriques in 1897 and proved by H.P.F.Swinnerton-Dyer in 1970 \cite{Sw5}.
\begin{thm}
A del Pezzo surface of degree 1, 5 or 7 must automatically have a $k$-point (See \textbf{\cite{Ct}} for a summary, and  \cite{Sko01} Corollary
3.1.5 for $d = 5$).
\end{thm}

\chapter{del Pezzo surfaces of degree 1}
\section{Rational points on del Pezzo surfaces of \\degree 1}
 Let $k=\bF_q$ be a finite field of characteristic $p$ and
let $k[x,y,z,w]$ be the weighted graded ring where the variables $x, y,
 z, w$ have weights 1, 1, 2, 3, respectively.

  Set $\mathbb P(1,1,2,3) :=$ Proj $k[x,y,z,w]$. Every del Pezzo
surface of degree 1 (called a \textbf{dP1}) over $k$ is isomorphic to a smooth sextic
in $\mathbb P(1, 1, 2, 3)$, and conversely \cite{Kollar1}. Thus, a
dP1 is given by an equation of the form
$$w^2 + H(x,y)w + z^3 + wz L(x,y) + Q(x,y)z^2+ G(x,y)z + F(x,y) =0,$$ where H, L, G and F are binary
homogeneous forms of degrees 3, 1, 4 and 6, respectively. Over fields of
char $k\neq 2$ or $3$, we may assume that $H(x,y)=Q(x,y)=0$.

Every dP1 comes endowed with a $\bF_q$-rational point--the base
point of the anticanonical linear system. A natural question arises,
``For what $q$ and what del Pezzo surfaces $X$ of degree 1, is it
true that $|X(\bF_q)|=1$?''

By Weil's theorem, if X has a unique $\bF_q$ point, i.e. $|X(\mathbb F_q)|=1=q^2+ qTrF^* +1$, the corresponding conjugacy class $\sigma \in
\mathscr{W}$(the Weyl group) must satisfy $ Tr \sigma =-q$. Analyzing the information given by Urabe's table \cite{Urabe} for $W(E_8)$,
we deduce that the unique rational point situation can
happen in the following cases:
\begin{itemize}
\item  Over $\mathbb F_2$, the conjugacy classes $C_{12}$, $C_{19}$ and $C_{48}$.
\item  Over $\mathbb F_3$, the conjugacy classes $C_{1}$, $C_{9}$, $C_{10}$ and $C_{38}$.
\item  Over $\mathbb F_4$, the conjugacy classes $C_{15}$.
\item  Over $\mathbb F_5$, the conjugacy classes $C_{44}$.
\item  Over $\mathbb F_7$, the conjugacy classes $C_{8}$.
\end{itemize}
The trace information shows that X may possibly have a unique rational point
over all the above five fields.

Geometric analysis can further reduce the number of possible situations.

  Blowing up the base point of $|{-K}_X|$, we obtain an elliptic
fibration map $\rho: X\dashrightarrow \bP^1$. The surface $X$ has a
unique $\bF_q$-rational point if and only if every fiber of $\rho$
contains exactly one $\bF_q$-point. Using the previous notation, the
fiber above $[m : n] \in \mathbb P^1$ is isomorphic to the elliptic
curve $w^2+z^3+G(t,1)z+F(t,1)=0$ where $t=\frac{m}{n}$. For smooth elliptic fibers, Hasse's theorem on elliptic fibers provides us the following estimation \cite{Manin}:
\begin{thm}(Hasse)
Let E be an elliptic curve defined over a finite field $\bF_q$ with
$q=p^r$ elements (p $\in \mathbb Z$ is a prime). The following
statement gives a bound of the size of $E(\bF_q)$,i.e. the number of points of E defined over $\bF_q$, which is denoted
by $N_q$.
$$N_q \ge q+1 - 2\sqrt{q}$$
\end{thm}

Hence, when $q\geq 5$, we have at least two $\bF_q$-points on each
smooth fiber.

But this elliptic fiber can be singular. Recall that a curve is said to be \textbf{singular} if there is a singular point on
the curve at which all the partial derivatives of the defining equation of the curve are zero. Therefore, singular elliptic fiber could have a cusp or a node.
\begin{lem}
   Suppose $C$ is a singular fiber,
   \begin{itemize}
   \item If $C$ has a cusp then the cusp point is of multiplicity bigger than 1. Hence it is defined over $\bF_q$.
   \item If $C$ has a node then the intersection point has multiplicity 2. It is also defined over $\bF_q$.
   \item Since the equation of the fiber corresponds to the equation of our dP1, we can preclude the case where $C$ consists of three incident lines.
   \item If $C$ consists of a line and a conic and they intersect, then both intersection points are defined over $\bF_q$ or the pair of intersection points is defined over $\bF_q$. If they do not intersect, then the line is defined over $\bF_q$, which gives us many $\bF_q$-rational points.
   \end{itemize}
\end{lem}
In each case, we can easily get more than one $\bF_q$-points on it. Therefore, in all the cases, $X$ will contain more than one rational points
on the fibers, hence on the surface.

 So $|X(\mathbb F_q)|=1$
can only happen when $q < 5$, in other words, over the fields $
\mathbb F_2$, $\mathbb F_3$ and $\mathbb F_4$.

A systematic, computer-aided search then allows us to determine all
dP1s with a unique $\mathbb F_q$-point. Our main result is:
\begin{thm}{\label{dP1}}
 Let X be a del Pezzo surface of degree 1 defined over $\mathbb F_q$. We have:
\begin{itemize}
       \item Over $\mathbb F_2$, there are two del Pezzo surfaces of degree 1 up to projective
       equivalence with a unique rational point, and the defining equations are:
             \begin{enumerate}
                \item $w^2+ z^3+w(x^3+x^2y+y^3)+z(x^4+x^3y+y^4)+(x^6+x^4y^2+y^6)=0$
                \item  $w^2+ z^3+w(x^3+x^2y+y^3)+z(x^4+x^2y^2+y^4)+(x^6+x^2y^4+y^6)=0$
             \end{enumerate}
       \item Over $\mathbb F_3$, there is a unique del Pezzo surfaces of degree 1 up to projective
       equivalence with a unique rational point, and its defining equation is:
               $$w^2+z^3+z(2x^4+x^2y^2+2y^4)+x^6+2x^4y^2+y^6=0.$$
       \item Over $\mathbb F_q$, where $q = 4$, del Pezzo surfaces of degree 1 have at least $q+1=5$ rational points. Furthermore, when $q \geq 5$, del Pezzo surfaces of degree 1 have at least $2q+1$ rational points.
\end{itemize}
\end{thm}
\section{The proof of Theorem \ref{dP1}}
\begin{proof}
Following the previous analysis, we only need to check the existence of the possible dP1s with a unique rational point
over three finite fields $\mathbb F_2$, $\mathbb F_3$ and $\mathbb F_4$. We
analyze them one by one. Start by analyzing dP1s over $\mathbb F_2$:
\begin{itemize}
\item \textbf{Find del Pezzo Surface of degree 1 over $\mathbb F_2$ with a
unique rational point}

 Over $\mathbb F_2$, del Pezzo surfaces of
degree 1 in  $\mathbb P(1, 1, 2, 3)$ are given by $$w^2 +
wzL(x,y)+wG_3(x,y)+z^3 + G_4(x, y)z + G_6(x, y) =0 $$ where $L(x,y)$
is a linear polynomial, and $G_3(x,y), G_4(x,y), G_6(x,y)$ are binary
homogeneous polynomials in x,y of degrees 3,4,6 respectively.

 But over $\mathbb F_2=\{ \bar 0,\bar 1\}$, we will have the relations: $ w^n=w, z^n=z, x^n=x$
and $y^n=y$ where $ n=1, \ldots, 6$, hence the equation becomes $$ w + z +
wz(px+qy)+w(ax+bxy+cy)+z(dx+exy+fy)+gx+hxy+iy=0$$

Analyzing of the coefficients g, h and i:
\begin{enumerate}
\item If $g=h=i=1$, then we will always have $[0,0,1,1]$ as a
rational point.
\item If one of g, h, i is 0, then
\begin{enumerate}
\item  $i=0\Longrightarrow [0,1,0,0]$ and [0,0,1,1] are distinct rational points.
\item   $g=0\Longrightarrow [1,0,0,0]$ and [0,0,1,1] are distinct rational
points.
\item   $h=0\Longrightarrow [1,1,0,0]$ and [0,0,1,1] are distinct rational
points.
\end{enumerate}
\item If two of g, h, i are 0, then
\begin{enumerate}
\item  $h=g=0$, $i=1\Longrightarrow [1,0,0,0]$ and [0,0,1,1] are distinct rational
points.
\item  $h=i=0$, $g=1\Longrightarrow [0,1,0,0]$ and [0,0,1,1] are distinct rational
points.
\item   $g=i=0$, $h=1\Longrightarrow [1,0,0,0]$ and [0,1,0,0] are distinct rational
points.
\end{enumerate}
\item If $g=h=i=0$ then [1,1,0,0] and [0,1,0,0] are distinct rational points.
\end{enumerate}
From the preceding analysis, $g=h=i=1$  are the only values that have a unique solution to our equation. In this case, we can write our
equation as
$$ w+z+x+xy+y+pwzx+qwzy+awx+bwxy+cwy+dzx+ezxy+fzy=0$$ with
a, b, c, e, d, f, p, q as coefficients, and x, y, z, w as variables.

The following pseudo-code finds all the solutions:

\emph{\textbf{Step 1} Set solution number =0, define all variables.\\
      \textbf{Step 2 }Repeat x, y, z, w from 0 to 1, hence run over all the possible
      pairs; repeat a, b, c, e, d, f, p, q from 0 to 1, hence run over all the
      possible combinations for coefficients.\\
      \textbf{Step 3} Define $func1 =w+z+x+xy+y+pwzx+qwzy+awx+bwxy+cwy+dzx+ezxy+fzy$. \\
      \textbf{Step 4} If func1 mod $2\equiv0$ then first count the current solution; then go back to Step
      2.
 } \\
\emph{\textbf{Step 5} If the number of solution is still 0, then
print the equation with no solution; otherwise, print the equation with a unique solution under linear equivalence. }

\textbf{Result} There exists a equation with a unique solution. Using
the above notation, the equation is $$
w+z+w(x+xy+y)+z(x+xy+y)+(x+xy+y)=0,$$
 with a unique solution [0,0,1,1].

Recall that a dP1 is given by $$w^2 +
wzL(x,y)+wG_3(x,y)+z^3 + G_4(x, y)z + G_6(x, y) =0. $$ There are 2 coefficients for $L(x,y)$, 4 coefficients for $G_3(x,y)$, 5 coefficients for $G_4(x,y)$ and 7 for $G_6(x, y)$. Therefore all the equations form an affine linear subspace of dimension $2+4+5+7=18$.

By the previous analysis, there is no mixed term as $wzL(x,y)$, which means $L(x,y)\equiv 0$. This reduces the dimension to 16. In order for the equation to have a unique solution, the following equations for $(x,y) = \{(0, 1), (1, 0), (1, 1)\}$ need to be satisfied:
 \begin{numcases}{}
 \nonumber
  G_3(x,y) \equiv 1, \\\nonumber
  G_4(x,y) \equiv 1, \\\nonumber
  G_6(x,y) \equiv 1.
 \end{numcases}\nonumber

These impose $3 \times 3=9$ conditions, hence the dimension is reduced to $ 16-9=7$. This gives us a list of $2^7=128$ equations of dP1s
with a unique $\mathbb F_2$ point.

Furthermore, we can fix the term $wG_3(x,y)$ to be $w(x^3+y^3+ x^2y)$, therefore we have 128/2=64 equations in the list.

\begin{claim}
 There are 32 equations in the above list that define smooth varieties, and among which, there are two non-isomorphic del
 Pezzo surfaces of degree 1 up to projective equivalence.
\end{claim}
\begin{proof}

First, let's check the smoothness of all the surfaces defined by the above 64 equations in the list. Let X be the
surface defined by the equation $$ w^2+ z^3+ w(x^3+y^3+ x^2y)+ z(x^4+ y^4+ x^3y)+ (x^6+x^5y+y^6).$$  The smoothness of X is
determined by the smoothness of X on each affine piece of $\mathbb
P(1,1,2,3)$.

For example, consider X in the fourth affine piece, which has a basis
$[a = \displaystyle \frac{x^3}{w}, b = \displaystyle \frac{x^2y}{w}, c = \displaystyle \frac{xy^2}{w}, d = \displaystyle \frac{y^3}{w}, e = \displaystyle \frac{xz}{w}, f = \displaystyle \frac{yz}{w}, g = \displaystyle \frac{z^3}{w^2}]$.
If the affine scheme given by translating X into this patch is nonsingular, then X is
smooth on this affine piece. If X are nonsingular on all four affine
patches, then X is smooth overall.

A simple calculation on X shows that X is a smooth surface, hence is a del Pezzo surface of degree 1.

Recall dP1 is given by $ w^2+ z^3+ w(x^3+y^3+ x^2y)+ zG_4(x, y) + G_6(x, y) =0. $
Checking each of those 64 equation in the same way, we get 32 equations defining smooth surfaces:
\newcounter{tmpcounter}
\begin{itemize}

\item When $G_4(x, y)$ contain 3 terms and $G_6(x, y)$ contain 5 or 7 terms , we have the following 17 equations:
    \begin{enumerate}

        \item $w^2 + z^3 + wy^3 + zy^4 + y^6 + zy^3x + wyx^2 + y^4x^2 + wx^3 + zx^4 +
    y^2x^4 + yx^5 + x^6$
        \item $w^2 + z^3 + wy^3 + zy^4 + y^6 + y^5x + wyx^2 + y^4x^2 + wx^3 + zyx^3 +
    zx^4 + yx^5 + x^6$
        \item $w^2 + z^3 + wy^3 + zy^4 + y^6 + wyx^2 + y^4x^2 + wx^3 + zyx^3 + y^3x^3
    + zx^4 + y^2x^4 + x^6$
        \item $w^2 + z^3 + wy^3 + zy^4 + y^6 + y^5x + wyx^2 + zy^2x^2 + wx^3
    + y^3x^3 + zx^4 + y^2x^4  + x^6$

        \item $w^2 + z^3 + wy^3 + zy^4 + y^6 + wyx^2 + y^4x^2 + wx^3 + zyx^3 + zx^4 +
    y^2x^4 + yx^5 + x^6$

        \item $w^2 + z^3 + wy^3 + zy^4 + y^6 + y^5x + wyx^2 + y^4x^2 + wx^3 + zyx^3 +
    y^3x^3 + zx^4 + x^6$

        \item $w^2 + z^3 + wy^3 + zy^4 + y^6 + wyx^2 + zy^2x^2 + y^4x^2 + wx^3 +
    y^3x^3 + zx^4 + yx^5 + x^6$

        \item $w^2 + z^3 + wy^3 + zy^4 + y^6 + wyx^2 + wx^3 + zyx^3 + y^3x^3 + zx^4 +
    y^2x^4 + yx^5 + x^6$

        \item $w^2 + z^3 + wy^3 + zy^4 + y^6 + zy^3x + y^5x + wyx^2 + wx^3 + y^3x^3 +
    zx^4 + yx^5 + x^6$

        \item $w^2 + z^3 + wy^3 + zy^4 + y^6 + zy^3x + wyx^2 + wx^3 + y^3x^3 + zx^4 +
    y^2x^4 + yx^5 + x^6$

        \item $w^2 + z^3 + wy^3 + zy^4 + y^6 + y^5x + wyx^2 + zy^2x^2 + wx^3 + y^3x^3
    + zx^4 + y^2x^4 + x^6$
        \item $w^2 + z^3 + wy^3 + zy^4 + y^6 + zy^3x + y^5x + wyx^2 + y^4x^2 + wx^3 +
    zx^4 + yx^5 + x^6$
        \item $w^2 + z^3 + wy^3 + zy^4 + y^6 + y^5x + wyx^2 + wx^3 + zyx^3 + y^3x^3 +
    zx^4 + yx^5 + x^6$
        \item $w^2 + z^3 + wy^3 + zy^4 + y^6 + y^5x + wyx^2 + zy^2x^2 + wx^3 + zx^4 +
    y^2x^4 + yx^5 + x^6$

        \item $w^2 + z^3 + wy^3 + zy^4 + y^6 + zy^3x + y^5x + wyx^2 + y^4x^2 + wx^3 +
    y^3x^3 + zx^4 + x^6$
        \item $w^2 + z^3 + wy^3 + zy^4 + y^6 + zy^3x + wyx^2 + y^4x^2 + wx^3 + y^3x^3
    + zx^4 + y^2x^4 + x^6$
        \item $w^2 + z^3 + wy^3 + zy^4 + y^6 + y^5x + wyx^2 + zy^2x^2 + y^4x^2 + wx^3
    + zx^4 + y^2x^4 + x^6$
        \setcounter{tmpcounter}{\theenumi}
    \end{enumerate}
 \item When $G_4(x, y)$ contain 5 terms and $G_6(x, y)$ contain 3, 5 or 7 terms , we have the following 8 equations:
    \begin{enumerate}
    \setcounter{enumi}{\thetmpcounter}
    \item $w^2 + z^3 + wy^3 + zy^4 + y^6 + zy^3x + y^5x + wyx^2 + zy^2x^2 +
    y^4x^2 + wx^3 + zyx^3 + y^3x^3 + zx^4 + x^6$
    \item $w^2 + z^3 + wy^3 + zy^4 + y^6 + zy^3x + y^5x + wyx^2 + zy^2x^2 + wx^3
    + zyx^3 + y^3x^3 + zx^4 + yx^5 + x^6$
    \item $w^2 + z^3 + wy^3 + zy^4 + y^6 + zy^3x + wyx^2 + zy^2x^2 + y^4x^2 +
    wx^3 + zyx^3 + zx^4 + y^2x^4 + yx^5 + x^6$
    \item $w^2 + z^3 + wy^3 + zy^4 + y^6 + zy^3x + wyx^2 + zy^2x^2 + wx^3 +
    zyx^3 + zx^4 + y^2x^4 + x^6$
    \item $w^2 + z^3 + wy^3 + zy^4 + y^6 + zy^3x + y^5x + wyx^2 + zy^2x^2 + wx^3
    + zyx^3 + zx^4 + x^6$
    \item $w^2 + z^3 + wy^3 + zy^4 + y^6 + zy^3x + wyx^2 + zy^2x^2 + y^4x^2 +
    wx^3 + zyx^3 + y^3x^3 + zx^4 + y^2x^4 + x^6$
    \item $w^2 + z^3 + wy^3 + zy^4 + y^6 + zy^3x + wyx^2 + zy^2x^2 + wx^3 +
    zyx^3 + y^3x^3 + zx^4 + y^2x^4 + yx^5 + x^6$
    \item $w^2 + z^3 + wy^3 + zy^4 + y^6 + zy^3x + y^5x + wyx^2 + zy^2x^2 +
    y^4x^2 + wx^3 + zyx^3 + zx^4 + yx^5 + x^6 $
    \setcounter{tmpcounter}{\theenumi}
    \end{enumerate}
 \item When $G_4(x, y)$ and $G_6(x, y)$ both contain 3 terms, we have the following 7 equations:
      \begin{enumerate}
       \setcounter{enumi}{\thetmpcounter}
        \item  $w^2+ z^3+w(x^3+x^2y+y^3)+z(x^4+x^3y+y^4)+(x^6+xy^5+y^6)=0$
        \item $w^2+ z^3+w(x^3+x^2y+y^3)+z(x^4+x^3y+y^4)+(x^6+x^4y^2+y^6)=0$
        \item  $w^2+ z^3+w(x^3+x^2y+y^3)+z(x^4+x^2y^2+y^4)+(x^6+x^2y^4+y^6)=0$
        \item  $w^2+ z^3+w(x^3+x^2y+y^3)+z(x^4+x^2y^2+y^4)+(x^6+x^3y^3+y^6)=0$
        \item $w^2+ z^3+w(x^3+x^2y+y^3)+z(x^4+x^2y^2+y^4)+(x^6+x^5y+y^6)=0$
        \item $w^2+ z^3+w(x^3+x^2y+y^3)+z(x^4+xy^3+y^4)+(x^6+xy^5+y^6)=0$
        \item $w^2+ z^3+w(x^3+x^2y+y^3)+z(x^4+xy^3+y^4)+(x^6+x^4y^2+y^6)=0$

      \end{enumerate}
 \end{itemize}

Next, let's check the equivalence between them. In order to
determine the isomorphisms between those equations, we need to list
all the transformations that fix $w^2 + z^3 + w(x^3+ y^3+ x^2y)$. They are as
follows:
\begin{equation}
\left.
(1)
\right.
\left\{
\begin{aligned}
x & \mapsto  y , \\
   y & \mapsto x+y , \\
 z & \mapsto mx^2+ny^2+z , \\    w & \mapsto w.
\end{aligned}
\right.
\left.
(2)
\right.
\left\{
\begin{aligned}
 x  & \mapsto x+y , \\
  y &  \mapsto  x, \\
 z & \mapsto mx^2+ny^2+z , \\    w & \mapsto w.
\end{aligned}
\right.
\nonumber
\end{equation}
\begin{equation}
\left.
(3)
\right.
\left\{
\begin{aligned}
 x  & \mapsto  y , \\
  y & \mapsto x+y , \\
  z & \mapsto  z, \\
  w & \mapsto p\cdot x^3+q \cdot y^3+r \cdot x^2y+ s \cdot xy^2 + w
\end{aligned}
\right.
\left.
(4)
\right.
\left\{
\begin{aligned}
 x  & \mapsto x+y , \\
  y &  \mapsto  x, \\
  z & \mapsto  z, \\
  w & \mapsto px^3+qy^3+r x^2y+ s xy^2 + w
\end{aligned}
\right.
\nonumber
\end{equation}
  Where $p, q, r, s$ are 0 or 1.

Applying the first two transformations to the above 32 equations, we get the following isomorphisms:
 \begin{enumerate}
    \item Equation \#1 is isomorphic to equation \#22 via transformation (1), isomorphic to \#2 via (2).
    \item Equation \#3 is isomorphic to \#15 via (1), isomorphic to \#24 via (2). And \#12 is isomorphic to \#24 via (1), isomorphic to \#26 via (2).
          \#20 is isomorphic to \#26 via (1), isomorphic to \#15 via (2).

          Therefore, \#3 is isomorphic to \#15, \#24, \#12, \#26 and \#20.
    \item  Equation \#4 is isomorphic to \#7 via (1), isomorphic to \#29 via (2). And \#11 isomorphic to \#29 via (2). Therefore \#4 is isomorphic to \#7, \#29 and \#11.
    \item  Equation \#5 is isomorphic to \#31 via (1), isomorphic to \#25 via (2).
    \item  Equation \#6 is isomorphic to \#10 via (1), isomorphic to \#23 via (2).
    \item  Equation \#8 is isomorphic to \#16 via (1), isomorphic to \#18 via (2).
    \item  Equation \#9 is isomorphic to \#19 via (1), isomorphic to \#13 via (2).
    \item  Equation \#14 is isomorphic to \#17 via (1), isomorphic to \#30 via (2).
    \item  Equation \#21 is isomorphic to \#27 via (1), isomorphic to \#32 via (2).
    \item  Equation \#28 is fixed under (1) and (2).
 \end{enumerate}
  Hence there are 10 dP1 with a unique rational point left after applying transformations (1) and (2) over $\mathbb F_2$:
  \begin{itemize}
    \item \#1: $w^2 + z^3 + wy^3 + zy^4 + y^6 + zy^3x + wyx^2 + y^4x^2 + wx^3 + zx^4 +
    y^2x^4 + yx^5 + x^6$

        \item \#3: $w^2 + z^3 + wy^3 + zy^4 + y^6 + wyx^2 + y^4x^2 + wx^3 + zyx^3 + y^3x^3
    + zx^4 + y^2x^4 + x^6$
        \item \#4: $w^2 + z^3 + wy^3 + zy^4 + y^6 + y^5x + wyx^2 + zy^2x^2 + y^4x^2 + wx^3
    + y^3x^3 + zx^4 + y^2x^4 + yx^5 + x^6$

        \item \#5: $w^2 + z^3 + wy^3 + zy^4 + y^6 + wyx^2 + y^4x^2 + wx^3 + zyx^3 + zx^4 +
    y^2x^4 + yx^5 + x^6$

        \item \#6: $w^2 + z^3 + wy^3 + zy^4 + y^6 + y^5x + wyx^2 + y^4x^2 + wx^3 + zyx^3 +
    y^3x^3 + zx^4 + x^6$

        \item \#8: $w^2 + z^3 + wy^3 + zy^4 + y^6 + wyx^2 + wx^3 + zyx^3 + y^3x^3 + zx^4 +
    y^2x^4 + yx^5 + x^6$

        \item \#9: $w^2 + z^3 + wy^3 + zy^4 + y^6 + zy^3x + y^5x + wyx^2 + wx^3 + y^3x^3 +
    zx^4 + yx^5 + x^6$
        \item \#14: $w^2 + z^3 + wy^3 + zy^4 + y^6 + y^5x + wyx^2 + zy^2x^2 + wx^3 + zx^4 +
    y^2x^4 + yx^5 + x^6$
        \item \#27: $w^2+ z^3+w(x^3+x^2y+y^3)+z(x^4+x^3y+y^4)+(x^6+x^4y^2+y^6)=0$
        \item  \#28: $w^2+ z^3+w(x^3+x^2y+y^3)+z(x^4+x^2y^2+y^4)+(x^6+x^2y^4+y^6)=0$
  \end{itemize}
 Next let's apply transformation (3) to these 10 equations:
 \begin{itemize}
 \item \#1 is isomorphic to \#18 via (3) when $p=r=s=0$, $q=1$, hence isomorphic to \#8.\\
         \#1 is isomorphic to \#19 via (3) when $p=s=0$, $q=r=1$, hence isomorphic to \#9.\\
         \#1 is isomorphic to \#23 via (3) when $q=r=p=0$, $s=1$, hence isomorphic to \#6.\\
          \#1 is isomorphic to \#24 via (3) when $p=s=q=1$, $r=0$, hence isomorphic to \#3.\\
          \#1 is isomorphic to \#21 via (3) when $p=r=0$, $q=s=1$, hence isomorphic to \#27.
 \item \#4 is isomorphic to \#14 via (3) when $p=q=r=0$, $s=1$.\\
     \#4 is isomorphic to \#28 via (3) when $s=q=r=0$, $p=1$.\\
 \end{itemize}
 \begin{lem}
 \#1 and \#4 are not isomorphic.
 \end{lem}
 \begin{proof}
 Applying all the linear transformations given above, we can easily get that those two equations are not isomorphic.
 \end{proof}
  In conclusion, there are two non-isomorphic dP1 up to projective equivalence, whose equations are given as \#1 and \#4.
\end{proof}

\item \textbf{Find del Pezzo Surface of degree 1 over $\mathbb F_3$ with a
unique rational point}

Let X denote a del Pezzo surface of degree 1 over $\mathbb F_3$,
then X is given by
$$ w^2 + z^3 + z^2G_2(x,y)+ G_4(x, y)z + G_6(x, y) =0,$$ where
$G_2(x,y), G_4(x,y)$, and $G_6(x,y)$ are binary homogeneous polynomials
in x, y of degree 2, 4 and 6 respectively. Over $\mathbb F_3$, we
have the relation $x^3=x, y^3=y, z^3=z$, where $x,y,z=\{ \bar 0,\bar
1,\bar 2\}$; hence we can write the equation explicitly as
follows:
$$ w^2+z+z^2(a x^2+bxy+ c y^2)+ z(dx^2+exy+fx^2y^2+gy^2)+(px^2+qxy+rx^2y^2+sy^2)=0 $$

The point [0, 0, 2, 2] is always a rational point on our surface over
$\mathbb F_3$. To avoid counting equivalent solutions
redundantly, we use the following approach to simplify the search.

Suppose X has a unique $\bF_3$-rational point [0,0,2,2]. Then the
equation of X at least shouldn't have [0, 1, z, w] as a solution, which is to say,
$$ w^2+z+z^2(c\cdot 1)+z(g\cdot 1)+s \neq 0. $$ Thus, we need to
\textbf{find  $(c, g, s)$ such that the equation $w^2+z+z^2c+zg+s=
0$ has no solution.} We use the following pseudo-code to establish
this: \\
\emph{\textbf{Step 1} Set solution number =0, define all variables.\\
      \textbf{Step 2 }Repeat w,z from 0 to 2, hence run over all the possible
      pairs; loop a,b,c from 0 to 2, or in other words, we run over all the
      possible combinations for coefficients.\\
      \textbf{Step 3} Define $func1 = w^2+ z+ az^2+ bz+ c$. \\
      \textbf{Step 4} If func1 mod $3\equiv 0$, then count the current solution and go back to Step
      2.
 } \\
\emph{\textbf{Step 5} If solution number is still 0, then print the
equation with no solution; otherwise, print the equation with a unique solution up to linear isomorphism.}

 \textbf{Result}: Only when $c=0$, $g=2$, $s=1$ the
equation has no solution. Plugging the value of these coefficient into the original equation, we get
$$w^2+z+z^2(ax^2+bxy)+z(dx^2+exy+fx^2y^2+2y^2)+(px^2+y^2+rx^2y^2+qxy)=0.$$
From the analysis above, we know that this is the only possible type of del Pezzo surfaces of degree 1 that can have a
unique rational point.

By means of varying over all possible values of other coefficients, we run a similar program as in the case
over $\bF_2$, to find solutions for this equation over $\mathbb
F_3$.

We have the following result:

\textbf{Result} The only coefficient combination that enables our equation to have a unique solution is $a = b = e= q =0, d=2, f=1, p=1, r=2$,
i.e., $$w^2+z+z(2x^2+ 0 \cdot xy+ x^2y^2+ 2y^2)+(x^2+y^2+2 \cdot x^2y^2+ 0 \cdot xy)=0.$$
However, $z(0 \cdot xy)$ can be $0$ , $z(x^3y + 2xy^3)$ or $z(2x^3y + xy^3)$.

 $2 \cdot x^2y^2$ can be $2 \cdot x^2y^4$, $2 \cdot x^4y^2$ or $\cdot x^2y^4+ \cdot x^4y^2$.

$0 \cdot xy$ can be $0$ , $x^5y + x^3y^3+ xy^5$, $x^5y + 2x^3y^3$, $x^5y  + 2xy^5$, $x^3y^3+ 2xy^5$, $2x^3y^3+ xy^5$, $2x^5y+ xy^5$, $2x^5y + x^3y^3$ or $2x^5y + 2x^3y^3+ 2xy^5$.

Therefore, we have $3 \times 3 \times 9 = 81$ equations.
\begin{lem}
Those 81 equations all define smooth del Pezzo surfaces of degree 1 over $\mathbb F_3$.
\end{lem}
\begin{proof}
We use the same smoothness checking program as in the case dP1 over $\mathbb F_2$ to check, but over field $\mathbb F_3$, then we can easily see they are all smooth.
\end{proof}

 \begin{lem}
 All del Pezzo surfaces of degree 1 over $\mathbb F_3$ have more than one rational point except for the following one:
$$w^2+z^3+z(2x^4+x^2y^2+2y^4)+x^6+2x^2y^4+y^6=0$$

 The surface is smooth, hence gives us a del Pezzo surface of degree 1.
\end{lem}
\begin{proof}
 First, we can fix $z(0 \cdot xy)$ to be $0$ or $z(2x^3y + xy^3)$, therefore we only have 54 equations left.

 \begin{itemize}
 \item  Let's look at the first 27 equations with $z(0 \cdot xy)= 0$:
 \begin{enumerate}
 \item $w^2 + z^3+ z(2x^4 + x^2y^2 + 2y^4) + (x^6 + y^6+ 2x^4y^2 )$
\item $w^2 + z^3+ z(2x^4 + x^2y^2 + 2y^4) + (x^6 + y^6+ 2x^4y^2 + x^5y + x^3y^3+ xy^5)$
\item $w^2 + z^3+ z(2x^4 + x^2y^2 + 2y^4) + (x^6 + y^6+ 2x^4y^2 + x^5y + 2x^3y^3)$
\item $w^2 + z^3+ z(2x^4 + x^2y^2 + 2y^4) + (x^6 + y^6+ 2x^4y^2 + x^5y  + 2xy^5)$
\item $w^2 + z^3+ z(2x^4 + x^2y^2 + 2y^4) + (x^6 + y^6+ 2x^4y^2 + x^3y^3+ 2xy^5)$
\item $w^2 + z^3+ z(2x^4 + x^2y^2 + 2y^4) + (x^6 + y^6+ 2x^4y^2 + 2x^3y^3+ xy^5)$
\item $w^2 + z^3+ z(2x^4 + x^2y^2 + 2y^4) + (x^6 + y^6+ 2x^4y^2 + 2x^5y  + xy^5)$
\item $w^2 + z^3+ z(2x^4 + x^2y^2 + 2y^4) + (x^6 + y^6+ 2x^4y^2 + 2x^5y + x^3y^3)$
\item $w^2 + z^3+ z(2x^4 + x^2y^2 + 2y^4) + (x^6 + y^6+ 2x^4y^2 + 2x^5y + 2x^3y^3+ 2xy^5)$

\item $w^2 + z^3+ z(2x^4 + x^2y^2 + 2y^4) + (x^6 + y^6+ 2x^2y^4 )$
\item $w^2 + z^3+ z(2x^4 + x^2y^2 + 2y^4) + (x^6 + y^6+ 2x^2y^4 + x^5y + x^3y^3+ xy^5)$
\item $w^2 + z^3+ z(2x^4 + x^2y^2 + 2y^4) + (x^6 + y^6+ 2x^2y^4 + x^5y + 2x^3y^3)$
\item $w^2 + z^3+ z(2x^4 + x^2y^2 + 2y^4) + (x^6 + y^6+ 2x^2y^4 + x^5y + 2xy^5)$
\item $w^2 + z^3+ z(2x^4 + x^2y^2 + 2y^4) + (x^6 + y^6+ 2x^2y^4 +  x^3y^3+ 2xy^5)$
\item $w^2 + z^3+ z(2x^4 + x^2y^2 + 2y^4) + (x^6 + y^6+ 2x^2y^4 + 2x^3y^3+ xy^5)$
\item $w^2 + z^3+ z(2x^4 + x^2y^2 + 2y^4) + (x^6 + y^6+ 2x^2y^4 + 2x^5y +xy^5)$
\item $w^2 + z^3+ z(2x^4 + x^2y^2 + 2y^4) + (x^6 + y^6+ 2x^2y^4 + 2x^5y +  x^3y^3)$
\item $w^2 + z^3+ z(2x^4 + x^2y^2 + 2y^4) + (x^6 + y^6+ 2x^2y^4 + 2x^5y + 2x^3y^3+ 2xy^5)$

\item $w^2 + z^3+ z(2x^4 + x^2y^2 + 2y^4) + (x^6 + y^6+ x^2y^4 + x^4y^2)$
\item $w^2 + z^3+ z(2x^4 + x^2y^2 + 2y^4) + (x^6 + y^6+ x^2y^4 + x^4y^2 + x^5y + x^3y^3+ xy^5)$
\item $w^2 + z^3+ z(2x^4 + x^2y^2 + 2y^4) + (x^6 + y^6+ x^2y^4 + x^4y^2 + x^5y + 2x^3y^3)$
\item $w^2 + z^3+ z(2x^4 + x^2y^2 + 2y^4) + (x^6 + y^6+ x^2y^4 + x^4y^2 + x^5y + 2xy^5)$
\item $w^2 + z^3+ z(2x^4 + x^2y^2 + 2y^4) + (x^6 + y^6+ x^2y^4 + x^4y^2 ++ x^3y^3+ 2xy^5)$

\item $w^2 + z^3+ z(2x^4 + x^2y^2 + 2y^4) + (x^6 + y^6+ x^2y^4 + x^4y^2  + 2x^3y^3+ xy^5)$
\item $w^2 + z^3+ z(2x^4 + x^2y^2 + 2y^4) + (x^6 + y^6+ x^2y^4 + x^4y^2 + 2x^5y + xy^5)$
\item $w^2 + z^3+ z(2x^4 + x^2y^2 + 2y^4) + (x^6 + y^6+ x^2y^4 + x^4y^2 + 2x^5y +  x^3y^3)$
\item $w^2 + z^3+ z(2x^4 + x^2y^2 + 2y^4) + (x^6 + y^6+ x^2y^4 + x^4y^2 + 2x^5y + 2x^3y^3+ 2xy^5)$
 \end{enumerate}

 Now let's find out the isomorphisms between them.
 \begin{lem}
    All the 27 equations above are isomorphic.
 \end{lem}
 \begin{proof}
  It's easy to get all the linear transformations between x, y, z and w which fix the term $w^2 + z^3+ z(2x^4 + x^2y^2 + 2y^4)$, and they are:
  \begin{equation}
\left.
(1)
\right.
\left\{
\begin{aligned}
x & \mapsto  x+2y , \\
   y & \mapsto x+y , \\
   z & \mapsto mx^2+ny^2+z , \\
   w & \mapsto w.
\end{aligned}
\right.
\left.
(2)
\right.
\left\{
\begin{aligned}
 x  & \mapsto x+2y , \\
  y &  \mapsto  2x+2y, \\
  z & \mapsto mx^2+ny^2+z , \\
   w & \mapsto w.
\end{aligned}
\right.
\left.
(3)
\right.
\left\{
\begin{aligned}
 x  & \mapsto x+y, \\
  y & \mapsto 2x+y, \\
z & \mapsto mx^2+ny^2+z , \\    w & \mapsto w.
\end{aligned}
\right.
\nonumber
\end{equation}
  \begin{equation}
\left.
(4)
\right.
\left\{
\begin{aligned}
 x  & \mapsto 2x+y, \\
  y &  \mapsto  x+y, \\
 z & \mapsto mx^2+ny^2+z , \\    w & \mapsto w.
\end{aligned}
\right.
\left.
(5)
\right.
\left\{
\begin{aligned}
x & \mapsto  2x+2y, \\
   y & \mapsto x+2y, \\
 z & \mapsto mx^2+ny^2+z , \\    w & \mapsto w.
\end{aligned}
\right.
\left.
(6)
\right.
\left\{
\begin{aligned}
 x  & \mapsto y, \\
  y &  \mapsto  x, \\
 z & \mapsto mx^2+ny^2+z , \\    w & \mapsto w.
\end{aligned}
\right.
\nonumber
\end{equation}
  \begin{equation}
\left.
(7)
\right.
\left\{
\begin{aligned}
 x  & \mapsto 2y, \\
  y & \mapsto x, \\
z & \mapsto mx^2+ny^2+z , \\    w & \mapsto w.
\end{aligned}
\right.
\left.
(8)
\right.
\left\{
\begin{aligned}
 x  & \mapsto  y, \\
  y &  \mapsto 2x, \\
 z & \mapsto mx^2+ny^2+z , \\    w & \mapsto w.
\end{aligned}
\right.
\left.
(9)
\right.
\left\{
\begin{aligned}
x & \mapsto   2y, \\
   y & \mapsto 2x, \\
 z & \mapsto mx^2+ny^2+z , \\    w & \mapsto w.
\end{aligned}
\right.
\nonumber
\end{equation}
  \begin{equation}
\left.
(10)
\right.
\left\{
\begin{aligned}
 x  & \mapsto x, \\
  y &  \mapsto 2y, \\
 z & \mapsto mx^2+ny^2+z , \\    w & \mapsto w.
\end{aligned}
\right.
\left.
(11)
\right.
\left\{
\begin{aligned}
 x  & \mapsto 2x, \\
  y & \mapsto y, \\
z & \mapsto mx^2+ny^2+z , \\    w & \mapsto w.
\end{aligned}
\right.
\left.
(12)
\right.
\left\{
\begin{aligned}
 x  & \mapsto  2x, \\
  y &  \mapsto 2y, \\
 z & \mapsto mx^2+ny^2+z , \\    w & \mapsto w.
\end{aligned}
\right.
\nonumber
\end{equation}
Where m, n are coefficients in $\mathbb F_3$.
Then we can check the isomorphisms between those 27 equations, and the result is:
\begin{itemize}
\item Via transformation (1),
 equation \#1 is isomorphic to \#9 when $m=0, n=2$,
  isomorphic to \#18 when $n=0, m=2$,
  isomorphic to \#27 when $m=n=1$.

      Via transformation (3), equation \#1 is isomorphic to \#2 when $m=0, n=2$, isomorphic to \#11 when $n=0, m=2$, isomorphic to \#20 when $m=n=1$.

      Via transformation (6), equation \#1 is isomorphic to \#10 when $m=0, n=0$.

\item Via transformation (3), equation \#3 is isomorphic to \#5 when $m=0, n=2$, isomorphic to \#15 when $n=0, m=2$, isomorphic to \#23 when $m=n=1$.

      Via transformation (4), equation \#3 is isomorphic to \#8 when $m=0, n=2$, isomorphic to \#19 when $n=0, m=2$, isomorphic to \#26 when $m=n=1$.

      Via transformation (5), equation \#3 is isomorphic to \#14 when $m=1, n=0$.

      Via transformation (6), equation \#3 is isomorphic to \#24 when $m=0, n=0$.

\item Via transformation (1), equation \#4 is isomorphic to \#6 when $m=0, n=2$, isomorphic to \#12 when $n=0, m=2$, isomorphic to \#21 when $m=n=1$.

      Via transformation (2), equation \#4 is isomorphic to \#9 when $m=0, n=0.$ $\Longrightarrow$ hence it is isomorphic to \#1.

      Via transformation (6), equation \#9 is isomorphic to \#16 when $m=0, n=0$.

      Via transformation (6), equation \#4 is isomorphic to \#14 when $m=0, n=0.$$\Longrightarrow$ hence it is isomorphic to \#2.

      Via transformation (8), equation \#4 is isomorphic to \#13 when $n=0, m=2$, isomorphic to \#22 when $m=n=1$. And \#22 is isomorphic to \#25 via exchange of x and y.

      Via transformation (10), equation \#4 is isomorphic to \#7 when $m=0, n=2$.

\item Via transformation (1), equation \#19 is isomorphic to \#1 when $m=1, n=0$,
\end{itemize}
Hence, they are all isomorphic to each other.

Therefore, we use $$w^2 + z^3+ z(2x^4 + x^2y^2 + 2y^4) + (x^6 + y^6+ 2x^4y^2 )$$ to represent the first 27 equation classes.
 \end{proof}

 \item Next, let's look at the last 27 equations with $z(0 \cdot xy)= z(2x^3y + xy^3)$.
 \begin{enumerate}
 \item $w^2 + z^3+ z(2x^4 + x^2y^2 + 2y^4 + 2x^3y + xy^3) + (x^6 + y^6+ 2x^4y^2 )$
\item $w^2 + z^3+ z(2x^4 + x^2y^2 + 2y^4 + 2x^3y + xy^3) + (x^6 + y^6+ 2x^4y^2 + x^5y + x^3y^3+ xy^5)$
\item $w^2 + z^3+ z(2x^4 + x^2y^2 + 2y^4 + 2x^3y + xy^3) + (x^6 + y^6+ 2x^4y^2 + x^5y + 2x^3y^3)$
\item $w^2 + z^3+ z(2x^4 + x^2y^2 + 2y^4 + 2x^3y + xy^3) + (x^6 + y^6+ 2x^4y^2 + x^5y+ 2xy^5)$
\item $w^2 + z^3+ z(2x^4 + x^2y^2 + 2y^4 + 2x^3y + xy^3) + (x^6 + y^6+ 2x^4y^2 + x^3y^3+ 2xy^5)$
\item $w^2 + z^3+ z(2x^4 + x^2y^2 + 2y^4 + 2x^3y + xy^3) + (x^6 + y^6+ 2x^4y^2 + 2x^3y^3+ xy^5)$
\item $w^2 + z^3+ z(2x^4 + x^2y^2 + 2y^4 + 2x^3y + xy^3) + (x^6 + y^6+ 2x^4y^2 + 2x^5y + xy^5)$
\item $w^2 + z^3+ z(2x^4 + x^2y^2 + 2y^4 + 2x^3y + xy^3) + (x^6 + y^6+ 2x^4y^2 + 2x^5y +  x^3y^3)$
\item $w^2 + z^3+ z(2x^4 + x^2y^2 + 2y^4 + 2x^3y + xy^3) + (x^6 + y^6+ 2x^4y^2 + 2x^5y + 2x^3y^3+ 2xy^5)$

\item $w^2 + z^3+ z(2x^4 + x^2y^2 + 2y^4 + 2x^3y + xy^3) + (x^6 + y^6+ 2x^2y^4 )$
\item $w^2 + z^3+ z(2x^4 + x^2y^2 + 2y^4 + 2x^3y + xy^3) + (x^6 + y^6+ 2x^2y^4 + x^5y + x^3y^3+ xy^5)$
\item $w^2 + z^3+ z(2x^4 + x^2y^2 + 2y^4 + 2x^3y + xy^3) + (x^6 + y^6+ 2x^2y^4 + x^5y + 2x^3y^3)$
\item $w^2 + z^3+ z(2x^4 + x^2y^2 + 2y^4 + 2x^3y + xy^3) + (x^6 + y^6+ 2x^2y^4 + x^5y + 2xy^5)$
\item $w^2 + z^3+ z(2x^4 + x^2y^2 + 2y^4 + 2x^3y + xy^3) + (x^6 + y^6+ 2x^2y^4  + x^3y^3+ 2xy^5)$
\item $w^2 + z^3+ z(2x^4 + x^2y^2 + 2y^4 + 2x^3y + xy^3) + (x^6 + y^6+ 2x^2y^4 + 2x^3y^3+ xy^5)$
\item $w^2 + z^3+ z(2x^4 + x^2y^2 + 2y^4 + 2x^3y + xy^3) + (x^6 + y^6+ 2x^2y^4 + 2x^5y + xy^5)$
\item $w^2 + z^3+ z(2x^4 + x^2y^2 + 2y^4 + 2x^3y + xy^3) + (x^6 + y^6+ 2x^2y^4 + 2x^5y +  x^3y^3)$
\item $w^2 + z^3+ z(2x^4 + x^2y^2 + 2y^4 + 2x^3y + xy^3) + (x^6 + y^6+ 2x^2y^4 + 2x^5y + 2x^3y^3+ 2xy^5)$

\item $w^2 + z^3+ z(2x^4 + x^2y^2 + 2y^4 + 2x^3y + xy^3) + (x^6 + y^6+ x^2y^4 + x^4y^2 )$
\item $w^2 + z^3+ z(2x^4 + x^2y^2 + 2y^4 + 2x^3y + xy^3) + (x^6 + y^6+ x^2y^4 + x^4y^2 + x^5y + x^3y^3+ xy^5)$
\item $w^2 + z^3+ z(2x^4 + x^2y^2 + 2y^4 + 2x^3y + xy^3) + (x^6 + y^6+ x^2y^4 + x^4y^2 + x^5y + 2x^3y^3)$
\item $w^2 + z^3+ z(2x^4 + x^2y^2 + 2y^4 + 2x^3y + xy^3) + (x^6 + y^6+ x^2y^4 + x^4y^2 + x^5y + 2xy^5)$
\item $w^2 + z^3+ z(2x^4 + x^2y^2 + 2y^4 + 2x^3y + xy^3) + (x^6 + y^6+ x^2y^4 + x^4y^2 + x^3y^3+ 2xy^5)$
\item $w^2 + z^3+ z(2x^4 + x^2y^2 + 2y^4 + 2x^3y + xy^3) + (x^6 + y^6+ x^2y^4 + x^4y^2 + 2x^3y^3+ xy^5)$
\item $w^2 + z^3+ z(2x^4 + x^2y^2 + 2y^4 + 2x^3y + xy^3) + (x^6 + y^6+ x^2y^4 + x^4y^2 + 2x^5y + xy^5)$
\item $w^2 + z^3+ z(2x^4 + x^2y^2 + 2y^4 + 2x^3y + xy^3) + (x^6 + y^6+ x^2y^4 + x^4y^2 + 2x^5y +  x^3y^3)$
\item $w^2 + z^3+ z(2x^4 + x^2y^2 + 2y^4 + 2x^3y + xy^3) + (x^6 + y^6+ x^2y^4 + x^4y^2 + 2x^5y + 2x^3y^3+ 2xy^5)$
 \end{enumerate}

  Now let's find out the isomorphisms between them.
 \begin{lem}
    All the 27 equations above are isomorphic.
 \end{lem}
 \begin{proof}
  It's easy to get all the linear transformations between x, y, z and w which fix the term $w^2 + z^3+ z(2x^4 + x^2y^2 + 2y^4 + 2x^3y + xy^3)$, and they are:
  \begin{equation}
\left.
(1)
\right.
\left\{
\begin{aligned}
x & \mapsto  2x+y , \\
   y & \mapsto x+2y , \\
   z & \mapsto mx^2+ny^2+z , \\
   w & \mapsto w.
\end{aligned}
\right.
\left.
(2)
\right.
\left\{
\begin{aligned}
 x  & \mapsto 2y , \\
  y &  \mapsto  x, \\
  z & \mapsto mx^2+ny^2+z , \\
   w & \mapsto w.
\end{aligned}
\right.
\nonumber
\end{equation}
 Where m, n are coefficients in $\mathbb F_3$.

Then we can check the isomorphisms between those 27 equations, and the result is:
\begin{itemize}
\item Via transformation (2), equation \#1 is isomorphic to \#2 when $m=1, n=2$, isomorphic to \#10 when $n=0, m=0$, isomorphic to \#21 when $m=1, n=0$.

      It is isomorphic to \#13 when $m=1, n=1$, isomorphic to \#5 when $n=1, m=0$, isomorphic to \#24 when $m=0, n=2$,

      isomorphic to \#25 when $m=2, n=1$, isomorphic to \#8 when $n=0, m=2$, isomorphic to \#18 when $m=n=2$.

\item Via transformation (2), equation \#19 is isomorphic to \#14 when $m=1, n=2$, isomorphic to \#3 when $m=1, n=0$.

      It is isomorphic to \#22 when $m=1, n=1$, isomorphic to \#3 when $n=1, m=0$, isomorphic to \#7 when $m=0, n=2$,

      isomorphic to \#25 when $m=2, n=1$ $\Longrightarrow$ hence isomorphic to \#1. It is also isomorphic to \#17 when $n=0, m=2$, isomorphic to \#9 when $m=n=2$.

\item   Via transformation (2), equation \#4 is isomorphic to \#2 when $m=1, n=2$ $\Longrightarrow$ hence isomorphic to \#1.

        It is also isomorphic to \#7 when $n=0, m=2$.

\item   Via transformation (2), equation \#12 is isomorphic to \#2 when $m=1, n=1$ $\Longrightarrow$ hence isomorphic to \#1.

        It is also isomorphic to \#16 when $n=1, m=2$.

\item   Via transformation (2), equation \#15 is isomorphic to \#2 when $m=2, n=2$ $\Longrightarrow$ hence isomorphic to \#1.

        It is also isomorphic to \#27 when $m=0, n=1$.

\item   Via transformation (2), equation \#20 is isomorphic to \#2 when $m=2, n=1$ $\Longrightarrow$ hence isomorphic to \#1.

\item   Via transformation (2), equation \#6 is isomorphic to \#3 when $m=2, n=0$ $\Longrightarrow$ hence isomorphic to \#1.

\item   Via transformation (2), equation \#23 is isomorphic to \#2 when $m=0, n=2$ $\Longrightarrow$ hence isomorphic to \#1.

\item   Via transformation (2), equation \#26 is isomorphic to \#24 when $m=n=0$ $\Longrightarrow$ hence isomorphic to \#1.

\end{itemize}
Hence, they are all isomorphic to each other. We can use the equation $$w^2 + z^3+ z(2x^4 + x^2y^2 + 2y^4 + 2x^3y + xy^3) + (x^6 + y^6+ x^2y^4 + x^4y^2 )$$ to represent this class.
 \end{proof}
\item At last, we can easily check that the representation equation of the first class $w^2 + z^3+ z(2x^4 + x^2y^2 + 2y^4) + (x^6 + y^6+ 2x^4y^2 )$ is isomorphic to the second one $w^2 + z^3+ z(2x^4 + x^2y^2 + 2y^4 + 2x^3y + xy^3) + (x^6 + y^6+ x^2y^4 + x^4y^2 )$ via the following transformation:
     \begin{equation}
        \left\{
    \begin{aligned}
    x & \mapsto x+ 2y , \\
   y & \mapsto  y , \\
   z & \mapsto  y^2+z , \\
   w & \mapsto w.
\end{aligned}
\right.
\nonumber
\end{equation}
 \end{itemize}
 Therefore, we have a unique dP1 over $\mathbb F_3$ with a unique rational point.
\end{proof}

\item \textbf{Find del Pezzo Surface of degree 1 over $\mathbb F_4$ with a
unique rational point}

Since $\mathbb F_4 = \displaystyle \frac{ \mathbb F_2 [t]}{t^2+t+1}$ is also of
characteristic 2, the equation is the same as over $\mathbb F_2$,
i.e. $$w^2 + wzL(x,y)+wG_3(x,y)+z^3 + G_4(x, y)z + G_6(x, y) =0 $$
Similar to the previous analysis, we always have [0,0,1,1] as the
base point.

Then we find all the equations that don't have [0,y,z,w] as their
solutions when $y\neq 0$. In other words, we have to find all \{
a,b,c,d \} over $\mathbb F_4$ such that the following equation has
no solution:
$$w^2+z^3+awyz+bwy^3+czy^4+dy^6=0$$
Running a similar computer program, we can find that the following equations satisfy the above requirement: \\
\begin{enumerate}
\item $w^2+z^3+\alpha wy^3+\alpha y^6=0$
\item $w^2+z^3+\alpha wy^3+ y^6=0$
\item $w^2+z^3+ wy^3+\alpha y^6=0$
\item $w^2+z^3+ wy^3+(1+\alpha)y^6=0$
\item $w^2+z^3+(1+\alpha) wy^3+  \alpha y^6=0$
\item $w^2+z^3+(1+\alpha) wy^3+(1+\alpha) y^6=0$
\end{enumerate}

Now, we need to compute all the solutions of these 6 equations over $\mathbb F_4$, and to see which corresponding original equation can have a unique rational point.

For example, let's consider an equation that contains $w^2+z^3+
wy^3+(1+\alpha) y^6=0$ as in case 4. The corresponding original equation for this
del Pezzo surface will be
$$ w^2+z^3+ awxz + bw+
cwx^2y+dwxy^2 + wy^3+ ezx+fzy + gzx^2y^2
$$ $$+kzxy^3+my^3
hx^2y+ixy^2+j+(1+\alpha)=0 $$

Then write $w=w_1+ \alpha w_2, x=x_1+ \alpha x_2,y=y_1+ \alpha
y_2,z=z_1+ \alpha z_2$. Use $\divideontimes$ to represent all
the coefficients and write $\divideontimes=\divideontimes_1+
\alpha\divideontimes_2 $, where all
$w_i,x_i,j_i,z_i,\divideontimes_i \in \mathbb F_2$.

Plug into the above original equation, then collect terms, in other words, to rewrite our equation in the form of $
F= F_1(x_i,y_i,z_i,w_i)+ \alpha  F_2(x_i,y_i,z_i,w_i)$. Then:

$F_1=w_1^2+w_2^2+z_1+z_2+z_1z_2+b_1w_1+b_2w_2+(a_1w_1+a_2w_2)(z_1x_1+z_2x_2)
 +(y_1+y_2+y_1y_2)(w_2+m_1+(k_1x_1+k_2x_2)z_1+(k_1x_2+k_2x_1+k_2x_2)z_2)
 +(a_1w_2+a_2w_1+a_2 w_2)(z_1 x_2+z_2 x_1+z_2
x_2)+(c_1w_1+c_2w_2)(x_1^2y_1+x_2^2y_1+x_2^2y_2)
 +(c_1w_2+c_2w_1+c_2w_2)(x_1^2y_2+x_2^2y_1)+(d_1w_1+d_2w_2)(x_1y_1^2+x_1y_2^2+x_2y_2^2)
 +(d_1w_2+d_2w_1+d_2w_2)(x_2y_1^2+x_1y_2^2) +(e_1x_1+e_2x_2)z_1
 +(e_1x_2+e_2x_1+e_2x_2)z_2
+(f_1y_1+f_2y_2)z_1+(f_1y_2+f_2y_1+f_2y_2)z_2
 +(z_1g_1+z_2g_2)((x_1^2+x_2^2)(y^2_1+y^2_2)+x^2_2+y^2_2)
 +(z_1g_2+z_2g_1+z_2g_2)(x^2_1y^2_2+x^2_2y^2_1+x^2_2y^2_2)
 +(h_1y_1+h_2y_2)(x^2_1+x^2_2)+(h_2y_1+h_1y_2+h_2y_2)x^2_2
 +(i_1x_1+i_2x_2)(y_1^2+y_2^2) +(i_2x_1+i_1x_2+i_2x_2)y_2^2+j_1; $

 $F_2
 =w_2^2+b_1w_2+b_2w_1+w_2b_2+(a_1w_1+a_2w_2)(z_1x_2+z_2x_1+z_2x_2)
 +(y_1+y_2+y_1y_2)(m_2+w_2+w_1+(k_1x_2+k_2x_2+k_2x_2)z_1+(k_1x_1+k_1x_2+k_2x_1)z_2)
 +(a_1 w_2+a_2 w_1+a_2 w_2) (z_1 x_2+z_2 x_1+z_1 x_1)
 +(c_1 w_2+c_2 w_2+c_2 w_1) (x_1 x_1 y_1+x_1 x_1 y_2+x_2 x_2 y_2)
 +(c_1 w_1+c_2 w_2) (x_1 x_1 y_2+x_2 x_2 y_1)
 +(d_1 w_1+d_2 w_2) (x_1 y_2 y_2+x_2 y_1 y_1)
 +(d_1 w_2+d_2 w_1+d_2 w_2) (x_2 y_1 y_1+x_1 y_1 y_1
 +x_2 y_2 y_2)+(e_1 x_2+e_2 x_1+e_2 x_2) z_1
 +(e_1 x_2+e_2 x_1+e_1 x_1) z_2+(f_1 y_2+f_2 y_1+f_2 y_2) z_1
 +(f_1 y_2+f_2 y_1+f_1 y_1) z_2
 +(x_1 y_1+x_2 y_2) (z_1 g_2+z_2 g_1+z_2 g_2)+(x_1 y_2+x_2 y_1+x_2
y_2)z_1 g_1+z_2 g_2)$$  +(h_1 y_2+h_2 y_1+h_2 y_2)(x_1^2+x_2^2)+(h_2
y_1+h_1 y_2+h_1 y_1)x_2 x_2+(i_1 x_2+i_2 x_1+i_2
x_2)(y_1^2+y_2^2)+(i_2 x_1+i_1 x_2+i_1 x_1)y_2^2+j_2+y_1+y_2+y_1
y_2$

Use an analysis similar to the one above to find the coefficients such that
$F_1=F_2=0$ has no common solution besides $[0,0,1,1]$.

It turns out we always have more than 2 solutions besides $[0,0,1,1]$  in case 4.

 Furthermore, we may use a very similar analysis and a computer program to
compute solutions of the original equations under the other 5 cases.

 \textbf{Result} All del Pezzo surfaces of degree 1 over $\mathbb F_4$
 have at least 2 rational points, hence at least 5
 rational points by Weil's theorem.

\item \textbf{del Pezzo Surface of degree 1 over $\mathbb F_5$ with
few rational points}

 From the previous geometric analysis using Hasse's theorem, we can easily conclude that over $\mathbb F_5$, if all the fibers are smooth
 elliptic curves, then we naturally have more than 1 point on
 the underlying surface. However, since there may be singular fibers
 of several different types, we still provide a computer proof as
 follows.

Over $\mathbb F_5$, it is given by $$ w^2 + z^3+ G_4(x, y)z + G_6(x,
y) =0$$and [0,0,1,2] is the base point as usual. Hence, we try to
find the equation that doesn't have any solution of the form [0,y,z,w] where $y\neq 0$. Which is to say, we try to find the
coefficients (e,s) such that
  $$ w^2 + z^3 + ey^4z+sy^6=0$$ doesn't have a solution.
Hence, we will find pairs (e,s) such that there is no nontrivial
solution to the above equation.

\emph{\textbf{Step 1} Set solution number =0, define all variables.\\
      \textbf{Step 2} Loop w,z from 0 to 5, y from 1 to 5, hence run over all the possible
      pairs; loop e,s from 0 to 5, hence run over all the
      possible combinations for coefficients.\\
      \textbf{Step 3} Define $func1 = w^2 + z^3 + ey^4z+sy^6$. \\
      \textbf{Step 4} If func1 mod $5\equiv 0$, then count current solution and go back to Step
      2.
 } \\
\emph{\textbf{Step 5} If the solution number is still 0, then print the
equation with no solution; otherwise, print the equation with a unique solution. } \\ \\
\textbf{Result:} The above equation always has solutions over any
pair of (e,s) over $\mathbb F_5$. Hence, all del Pezzo surfaces of
degree 1 over $\mathbb F_5$ will have at least 2 rational points.

Now, let's find a lower bound on the number of rational points on
dP1 over $\bF_5$. By Weil's theorem, we can get $N(\bF_q)\equiv
1 $(mod q) on the del Pezzo surface, hence we naturally have a lower bound as $q+1$. Can we improve the lower bound of the number of rational
solutions?

  Let's analyze the rational points on the fibers. By Hasse's theorem, when $q \geq 5$, we easily get that $N_q(\bF_q)\geq 2$ on each smooth fiber, and a singular fiber will have more rational points than a smooth fiber. Therefore, we may assume that all the fibers are smooth when computing the lower bound of rational points.

If all fibers over X are smooth, at least we have $(q+1) \cdot 1 +1=q+2$ points on the surface. By the Chevalley-Warning theorem or the Weil theorem, $N(\bF_q)\equiv
1 $(mod q). Therefore, $N(\bF_q)\geq 2q+1$ when $q \geq 5$.

\item \textbf{del Pezzo Surface of degree 1 over $\mathbb F_q$ where $q \geq 7$ with
few rational point}
\begin{claim}
 Over $\bF_q$ where $q \geq 7$, a del Pezzo surface of degree 1 has at least $2q+1$ rational
 points.
 \end{claim}
 \begin{proof}
 Similarly to the case over $\bF_5$, we get this result by analyzing the rational points on fibers. Since when $q \geq 7$, $N_q(\bF_q)\geq (1+q-2\sqrt{q})> 2$, we have at least 3 points on each fiber, then we have $(q+1) \cdot 2 +1=2q+3$ points on the surface. Therefore, $N(\bF_q)\geq 3q+1$ when $q \geq 7$.
 \end{proof}

\end{itemize}
\end{proof}

\chapter{del Pezzo surfaces of degree 2}
Let X be a del Pezzo surface of degree 2. Then X can be realized as
the surface in weighted projective space $\mathbb P(1, 1, 1, 2)$
given by the equation
\begin{equation}
w^2 + wG_2(x, y, z) + G_4(x, y, z)=0,
\end{equation}
where $G_2(x, y, z)$ and $G_4(x,y,z)$ are homogeneous polynomials of degree 2 and 4 respectively. When we consider X over a field of characteristic $\neq 2$, we may assume $G_2(x, y, z)=0$. So X is
a double cover of $\mathbb P^2$ ramified over the (smooth) quartic
curve $F = 0$, which is a canonical model of a curve of genus 3. The
lifts of the 28 bitangents to this quartic curve come in pairs,
forming 56 exceptional curves on X.

\section{Rational points on del Pezzo surfaces of \\ degree 2}
By Weil's theorem, and the information on traces from the Urabe's
table \cite{Urabe} for del Pezzo surface of degree 2, the unique
k-rational point situation can only happen in the following cases:
\begin{itemize}
\item  Over $\mathbb F_2$, Carter symbol: $5A_1$, with orbit decomposition $2^{4}\cdot 2^{8}\cdot 2^{16}$, $H^1\cong \mathbb Z_2 \times \mathbb Z _2$ \\
                        and Carter symbol: $A_3+3A_1$, with orbit decomposition $2^{2}\cdot 2^{6}\cdot 4^{10}$, $H^1\cong \mathbb Z_2 \times \mathbb Z_2$ \\
                        and Carter symbol: $2A_3+A_1$, with orbit decomposition $2^{4}\cdot 4^{12 }$, $H^1\cong \mathbb Z_4 \times \mathbb Z_4$.
\item  Over $\mathbb F_3$, Carter symbol: $D_4+3A_1$, with orbit decomposition $2^{10}\cdot 6^{6}$, $H^1\cong {\mathbb Z_2} ^4$.
\item  Over $\mathbb F_4$, Carter symbol: $6A_1$, with orbit decomposition $2^{12}\cdot 2^{16}$, $H^1\cong {\mathbb Z_2} ^4$.
\end{itemize}
\begin{thm} Let $X$ be a del Pezzo surface of degree 2 defined over $\mathbb F_q$. We have:
\begin{itemize}
       \item Over $\mathbb F_2$, there is no del Pezzo surface of degree 2 with a unique rational
       point.
       \item Over $\mathbb F_3$, there is a unique del Pezzo surface of degree 2 up to projective
       equivalence with a unique rational point, and its defining equation is:
        $$w^2 + z^4 + 2z^2y^2 + y^4 + z^3x + 2zx^3=0.$$
       \item Over $\mathbb F_q$, where $q \geq 4$, del Pezzo surfaces of degree 2 have at least $2q+1$ rational
       points.
\end{itemize}
\end{thm}

\section{The proof of Theorem 4.1}
\begin{proof}
 From the analysis above, we only need to check the possible existence of a dP2 with a unique rational
 point in the case of $\mathbb F_2$, $\mathbb F_3$, $\mathbb F_4$.

\begin{itemize}
 \item

 Consider the case where our del Pezzo surfaces of degree 2(\textbf{dP2}) are defined in $\mathbb P(1,1,1,2)$ over $\mathbb
F_2$ and $\mathbb F_4$. Here the equation of dP2 is $$ w^2+w
G_2(x,y,z)+ G_4(x,y,z)=0.$$
 Let $X(\mathbb F_q)= \{ (\mathbb F_q)$-solutions of the equation $w^2+w
G_2(x,y,z)+ G_4(x,y,z)=0 \} $ denote the rational points on our dP2 X, and
let $C= \{ L_2(x,y,z)=0 \}$ denote the zero locus of $L_2$.

We have $ \# X(\mathbb F_q) \geq \# C(\mathbb F_q)$, where $q=2^r$. In
order to find all the dP2 that contain a unique $\mathbb F_q$-point, we must
have $ \# C(\mathbb F_q)=0$ or 1. That is, $G_2(x,y,z)$ must be a
nonsplit quadric.

Over $\mathbb F_2$, without loss of generality, we may assume $G_2(x,y,z)=
x^2+ xy+ y^2$. We can see that $G_4(x,y,z)\equiv 1$ except for $x=y=z=0$.

 \textbf{Result}

 Running the similar computer program as in the analysis of dP1, we find 256
equations with a unique solution [0,0,1,0] over $\mathbb F_2$.

Since a dP2 is given by $w^2+w (x^2+ xy+ y^2) + G_4(x,y,z)=0  $, the equations of all dP2's form an affine subspace of dimension 15. However, the uniqueness of the solution requires that $G_4(x,y,z)\equiv 1$ except for $x=y=z=0$, which imposes 7 conditions on the coefficients of $G_4(x,y,z)$. Therefore, they form an affine linear subspace of dimension $15-7 =8$, which also gives us the list of $2^8$ equations. To save space, we omit the list here, since we are not going to use those according to the following lemma.

\begin{lem}
All the 256 equations above are not smooth. Hence there is no dP2 over
$\bF_2$ with a unique solution.
\end{lem}
\begin{proof}
We check the smoothness of each equation by checking it on all the
affine patches of $\bP(1,1,1,2)$.
\end{proof}

\item

 Consider the case over $\mathbb F_4$. Similar to the analysis over
$\mathbb F_2$, $G_2(x,y,z)$ must be a non-split quadric, and we
may assume $G_2(x,y,z)= x^2+ \alpha xy+ y^2$. To simplify the
procedure, we start to find all the equations that don't have
[0, y, z, w] as a solution. That is to say, $$w^2+y^2w+ G_4(0,y,z)=0$$
has no non-trivial solution over $\mathbb F_4$. Explicitly, we have
$$w^2+y^2w+ s_1 y^4+s_2 z^4+s_3 y^2z^2+s_4yz=0.$$

 Write $$s_i=a_i+\alpha b_i$$ $$w=w_1+\alpha w_2$$ $$z=z_1+\alpha z_2$$ $$y=y_1+\alpha y_2.$$
Plugging in and collecting terms, we rewrite the equation in form of
$F_1(y_i,z_i,w_i)+\alpha F_2(y_i,z_i,w_i)=0$. We have
$$ F_1 = w_1^2+ w_2^2+ w_1+a_1+a_2z_1+b_2z_2+a_3(z_1^2+z_2^2)+b_3z_2$$
 $$F_2 = w_1^2+w_2+b_1+b_2z_2+b_2z_1+a_2z_2+a_3z_2^2+b_3z_1^2+b_4z_1+b_4z_2+a_4z_2$$
 Run a similar computer program to find all the combination of coefficients such that
 $F_1=F_2=0$ have no common solution.

  \textbf{Result}

  There is no such combination. Hence, all dP2
surfaces have at least 2 rational points over $\mathbb F_4$, hence
they have more than 5 rational points.

\item

Consider the case over $\mathbb F_3$, where our equation
can be written as $w^2+G_4(x,y,z)=0$. Over $\mathbb F_3$, we have
relation $x^3=x$ for $x=\{ \bar 0, \bar 1, \bar 2\}$. Therefore, we can
write our equation more explicitly as
$$w^2+(ax^2+bxy+cxz)+dx^2y^2+ex^2yz+fx^2z^2+gxy^2z+(py^2+qyz+rz^2y^2+sz^2).$$

WLOG, we can assume the unique rational point is [0, 1, 0, 0] in $\mathbb P(1, 1, 1, 2)$.
Running a computer-aided program as before, the only combination of the coefficients which enables our equation to have a unique solution is
$$a=1, f=2, s=1, b=c=d=e=g=p=q=r=0$$
which is $ w^2 + x^4 + 0 \cdot xy + 0 \cdot xz + 2x^2z^2 + 0 \cdot yz + z^4 = 0. $
However, $0 \cdot xy$ can be $0\cdot x^3y + 0\cdot xy^3$, $x^3y + 2\cdot xy^3$ or $2\cdot x^3y + xy^3$.

We are in the same situation with $0 \cdot xz$ and $0 \cdot yz$, so we have 27 equations with a unique rational point [0, 1, 0, 0].

  Running the smoothness checking program similarly as before, we have the following 22 equations left which define smooth dP2s:
\begin{enumerate}
    \item $w^2 + z^4 + z^3y + 2zy^3 + y^3x + 2z^2x^2 + 2yx^3 + x^4 $
    \item $w^2 + z^4 + 2z^3y + zy^3 + y^3x + 2z^2x^2 + 2yx^3 + x^4 $
\item $w^2 + z^4 + 2z^3y + zy^3 + 2z^3x + 2z^2x^2 + zx^3 + x^4 $
\item $w^2 + z^4 + 2z^3y + zy^3 + z^3x + 2z^2x^2 + 2zx^3 + x^4 $
\item $w^2 + z^4 + 2z^3y + zy^3 + z^3x + 2y^3x + 2z^2x^2 + 2zx^3 + yx^3 +
    x^4 $
    \item $w^2 + z^4 + 2z^3y + zy^3 + z^3x + y^3x + 2z^2x^2 + 2zx^3 + 2yx^3 +
    x^4 $
\item $w^2 + z^4 + 2z^3x + 2y^3x + 2z^2x^2 + zx^3 + yx^3 + x^4 $
\item $w^2 + z^4 + 2z^3y + zy^3 + 2z^3x + 2y^3x + 2z^2x^2 + zx^3 + yx^3 +
    x^4 $

\item $w^2 + z^4 + 2y^3x + 2z^2x^2 + yx^3 + x^4 $
\item $w^2 + z^4 + 2z^3y + zy^3 + 2z^3x + y^3x + 2z^2x^2 + zx^3 + 2yx^3 +
    x^4 $
\item $w^2 + z^4 + z^3y + 2zy^3 + 2y^3x + 2z^2x^2 + yx^3 + x^4 $
\item $w^2 + z^4 + 2z^3y + zy^3 + 2y^3x + 2z^2x^2 + yx^3 + x^4 $
\item $w^2 + z^4 + z^3x + y^3x + 2z^2x^2 + 2zx^3 + 2yx^3 + x^4 $
\item $w^2 + z^4 + z^3y + 2zy^3 + z^3x + 2y^3x + 2z^2x^2 + 2zx^3 + yx^3 +
    x^4 $

    \item $w^2 + z^4 + z^3y + 2zy^3 + z^3x + y^3x + 2z^2x^2 + 2zx^3 + 2yx^3 +
    x^4 $

    \item $w^2 + z^4 + 2z^3x + y^3x + 2z^2x^2 + zx^3 + 2yx^3 + x^4 $
\item $w^2 + z^4 + z^3y + 2zy^3 + 2z^3x + 2y^3x + 2z^2x^2 + zx^3 + yx^3 +
    x^4$
    \item $w^2 + z^4 + z^3y + 2zy^3 + 2z^3x + y^3x + 2z^2x^2 + zx^3 + 2yx^3 +
    x^4 $

\item $w^2 + z^4 + z^3y + 2zy^3 + 2z^3x + 2z^2x^2 + zx^3 + x^4 $

    \item $w^2 + z^4 + z^3y + 2zy^3 + z^3x + 2z^2x^2 + 2zx^3 + x^4 $
\item $w^2 + z^4 + y^3x + 2z^2x^2 + 2yx^3 + x^4 $

\item $w^2 + z^4 + z^3x + 2y^3x + 2z^2x^2 + 2zx^3 + yx^3 + x^4 $
\end{enumerate}

Eliminate those equations from above list that are isomorphic to others via transformations $x \mapsto -x$, $y \mapsto -y$ and $z \mapsto -z$.

After this elimination, we have 6 equations left:
\begin{enumerate}
\item \#1 : $w^2 + z^4 + z^3y + 2zy^3 + y^3x + 2z^2x^2 + 2yx^3 + x^4 $
\item \#3 : $w^2 + z^4 + 2z^3y + zy^3 + 2z^3x + 2z^2x^2 + zx^3 + x^4 $
\item \#5 : $w^2 + z^4 + 2z^3y + zy^3 + z^3x + 2y^3x + 2z^2x^2 + 2zx^3 + yx^3 +
    x^4 $
\item \#6 : $w^2 + z^4 + 2z^3y + zy^3 + z^3x + y^3x + 2z^2x^2 + 2zx^3 + 2yx^3 +
    x^4 $
\item \#7 : $w^2 + z^4 + 2z^3x + 2y^3x + 2z^2x^2 + zx^3 + yx^3 + x^4 $
\item \#9 : $w^2 + z^4 + 2y^3x + 2z^2x^2 + yx^3 + x^4 $
\end{enumerate}
\begin{claim}
There is a unique dP2 over $\bF_3$ with a unique rational point.
\end{claim}
\begin{proof}
Considering the isomorphisms between the above 6 equations, we have:

$(3)\mapsto(1)$ by $x\mapsto z-x, z\mapsto x+z;$ \\
$(2)\mapsto(5)$ by $x\mapsto x-z, z\mapsto x+z; $\\
$(3)\mapsto(5)$ by $x\mapsto x-z, z\mapsto x+z, y\mapsto x+y+z$. \\
$(4)\mapsto(6)$ by $x\mapsto z-x, z\mapsto -x-z;$  \\
$(6)\mapsto(2)$ by $x\mapsto z-x, z\mapsto -x-z;$  \\
They are isomorphic to each other, hence we have a unique dP2 surface with a unique rational point over
$\bF_3$.
\end{proof}

\item Next we consider del Pezzo surfaces of degree 2 over $\mathbb F_q$ where $q \geq 4$ with
few rational points.
\begin{claim}
 Over $\bF_q$ where $q \geq 4$, a del Pezzo surface of degree 2 has at least $2q+1$ rational
 points.
 \end{claim}
 \begin{proof}

   For $q=4$, according to the information of trace in Urabe's table, the minimal trace of all the conjugacy classes of dP2 is -4. And the corresponding class is $C_2$ of Carter symbol: $6A_1$.
   Similar to the program for dP1 surface over $\bF_4$, we run a program to check solutions for this dP2 over $\bF_4$.

   The result shows that it has at least 6 solutions up to projective equivalence. By Weil's theorem, the number of solutions must equal to 1 modulo q. Hence a dP2 has at least $2q+1=9$ solutions over $\bF_4$.\\
   For $q>5$, we have $q+ Tr F^* > 1$ for all the other classes. Then by Weil's theorem again, $N(\bF_q)= q^2 + q  Tr F^*  + 1$, so we always have $N(\bF_q) \geq 2q+1$.

 \end{proof}

\end{itemize}

\end{proof}

\chapter{Rationality of cubic surfaces }

Let X be a nonsingular cubic surface defined
over a given field $k$. Hence X is a smooth cubic hypersurface $X
\subset \mathbb P^3$.

 It is well known that over the complex numbers any such X is birationally equivalent
to a projective plane. The problem of finding necessary and
sufficient conditions for X to be birationally equivalent to a
projective plane over $k$ was first raised by B.Segre\cite{Manin}. In the situation where
X is defined over an algebraic number field $k$, a complete answer
was given by H.P.F. Swinnerton-Dyer\cite{Sw3} in his paper ``The birationality
of cubic surfaces over a given field''.

 The geometry of X has been known and studied for a long time; X is isomorphic to the projective plane with six points in
general position blown up and there are 27 lines defined over the algebraic closure of
the ground field.

The configuration of the 27 lines gives implicit yet important
information about the arithmetic and geometry of the cubic surface.
Here is a natural question:
\begin{problem}
  Given a smooth cubic surface, can one find the equations of
the 27 lines, and the configuration of these lines?
\end{problem}
This is quite difficult, since usually the 27 lines are not even
defined over the field of definition of the surface. Throughout this chapter, we are interested in the geometry of the cubic surface over a
finite field of small cardinality.


\section{Eckardt points on cubic surfaces}

In the previous discussion of Corti-Koll\'ar-Smith example[KSC04], we have seen that there exists a unique cubic surface with a unique rational point over $\mathbb F_2$. Let's compute all the 27 lines on it, and the configuration of those lines.

 \begin{exam}

 Let X denote this unique cubic surface over $\bF_2$. Then all the 27 lines on $\bar X$ are NOT defined over the ground field $\bF_2$, even the unique Eckardt point is defined over $\bF_2$.

\begin{defn}
 An \textbf{Eckardt Point} of a non-singular cubic surface is a point through which pass
three coplanar lines of the surface.
\end{defn}

 By Weil's theorem, X corresponds to the conjugacy class of trace -2, in order to have a unique $\bF_2$-point. Therefore, X corresponds to the class $C_{11}$, with index 0, order 3 and the first cohomology group $\mathbb Z_3 \times \mathbb Z_3$. The minimum field extension of $\bF_2$ where all the 27 lines are defined is the splitting field of degree 3. This information is easy to get from Manin's table for cubic surface. Now let's compute the 27 lines on X over $\bF_8$.

 Let $\alpha$ denote a generator of the multiplicative group $\bF_8^{\times}$, and we can check that $\alpha$ also satisfies $\alpha^3+\alpha+1=0$ over $\bF_8$.
 Recall that the equation of X in $\bP^3$ is $xw^2 + x^2w + y^3 + z^2y + z^3 + x^3 + yzx + (y+z)x^2=0$.

It is hard to compute the equation of the 27 lines directly out of the equation of X. Instead we use Pl\"ucker coordinate to denote all the lines. By computing the Gr\"obner basis of the generated ideal, we get the following 27 lines:

\vspace{1cm}

\begin{tabular}{ll}
  \hline
  Number of the line & $[p_{01}, p_{02}, p_{03}, p_{23}, p_{31}, p_{12}]$ \\
  \hline
  1 & $[\alpha^6, \alpha^5, 0, 0, 0, 1]$ \\
  2 & $[\alpha^5, \alpha^3, 0, 0, 0, 1]$ \\
  3 & $[\alpha^3, \alpha^6, 0, 0, 0, 1]$ \\
  4 & $[\alpha^4, \alpha^3, \alpha^3,0,1,1]$ \\
  5 & $[\alpha^2, \alpha^5, \alpha^5,0,1,1]$ \\
  6 & $[\alpha, \alpha^6, \alpha^6, 0,1,1]$ \\
  7 & $[\alpha^6, \alpha^5, \alpha, \alpha^4, \alpha^3, 1]$ \\
  8 & $[\alpha^6, \alpha^4, \alpha^6, 1, 0, 1]$ \\
  9 & $[\alpha^5, \alpha^3, \alpha^2, 1, \alpha^6, \alpha]$ \\
  10 & $[\alpha^5, \alpha, \alpha^5, 1, 0, 1]$ \\
  11 & $[\alpha^4, \alpha^3, \alpha^4, 1, \alpha^4, \alpha]$ \\
  12 & $[\alpha^4, \alpha^2,\alpha, 1, 1, 1]$ \\
  \hline
\end{tabular}

\begin{tabular}{ll}
  \hline
  Number of the line & $[p_{01}, p_{02}, p_{03}, p_{23}, p_{31}, p_{12}]$ \\
  \hline
  13 & $[\alpha^4, \alpha, \alpha, 1, \alpha, \alpha^2]$ \\
  14 & $[\alpha^4, \alpha, 0, 1, \alpha^4, \alpha]$ \\
  15 & $[\alpha^3, \alpha^6, \alpha^4, 1, \alpha^5, \alpha^2]$ \\
  16 & $[\alpha^3, \alpha^2, \alpha^3, 1, 0, 1]$ \\
  17 & $[\alpha^2, \alpha^5, \alpha^2, 1, \alpha^2, \alpha^4]$ \\
  18 & $[\alpha^2, \alpha, \alpha^4, 1, 1, 1]$ \\
  19 & $[\alpha^2, 1, 1, 1, \alpha^4, \alpha]$ \\
  20 & $[\alpha^2, 1, 0, 1, \alpha^2, \alpha^4]$ \\
  21 & $[\alpha, \alpha^6, \alpha, 1, \alpha, \alpha^2]$ \\
  22 & $[\alpha, \alpha^4, \alpha^2, 1, 1, 1]$ \\
  23 & $[\alpha, 1, 1, 1, \alpha^2, \alpha^4]$ \\
  24 & $[\alpha, 1, 0, 1, \alpha, \alpha^2]$ \\
  25 & $[0, 0, 0, 0, 1, \alpha^4]$ \\
  26 & $[0, 0, 0, 0, 1, \alpha^2]$ \\
  27 & $[0, 0, 0, 0, 1, \alpha]$ \\
  \hline
\end{tabular}

 \vspace{1cm}
Now that we have all the 27 lines on X, we need one more step to make clear the configuration of these lines, in other words, we need to check whether they are skew or coplanar.

If $p_{ij}$  and $ p'_{ij}$ are the Pl\"ucker coordinates of two lines, then they are coplanar precisely when the following equation is satisfied:
    $$0 = p_{01}p'_{23} + p_{02}p'_{31} + p_{03}p'_{12} + p_{23}p'_{01} + p_{31}p'_{02} + p_{12}p'_{03}$$
If this equation is not satisfied, then they are skew.

In Pl\"ucker coordinates $[p_{01}, p_{02}, p_{03}, p_{23}, p_{31}, p_{12}]$, let \textbf{d} denote $[p_{01}, p_{02}, p_{03}]$ and \textbf{m} denote $[p_{23}, p_{31}, p_{12}]$.
Then two coplanar lines, neither of which contains the origin, have a common point
\begin{equation}
(x_0 : \textbf{x}) = (\textbf{d} \cdot \textbf{m¡ä}: \textbf{m} \times \textbf{m¡ä}).\nonumber
\end{equation}     	
Look at all the intersection points of these 27 lines over $\mathbb F_8$; there must be some Eckardt points. A simple calculation shows:
$$ \# X(\mathbb F_8)=8^2+7 \ 8+121.$$ Namely there are 135 points on pairs of lines. Suppose we have $N$ Eckardt points on the surface X, which satisfies:
\begin{eqnarray}
 \nonumber 135-2N &\leq& 121 \\\nonumber
 N &\geq& 9
\end{eqnarray}
After checking the two conditions for the 27 lines, we have the following:
\begin{itemize}
\item $l_1, l_2$ and $l_3$ are coplanar to each other, and provide an Eckardt point.
\item $l_{25}, l_{26}$ and $l_{27}$ are coplanar, and provide an Eckardt point.
 The first two cases give us the same point at infinity.

\item $l_4, l_5$ and $l_6$ are coplanar to each other, and provide an Eckardt point $[1: 0: 0: 0]$.
\item $l_{17}, l_{20}$ and $l_{23}$ are coplanar, and provide an Eckardt point$[1: 0: 0: 0]$.
\item $l_8, l_{10}$ and $l_{16}$ are coplanar, and provide an Eckardt point $[1: 0: 0: 0]$.
\item $l_{12}, l_{18}$ and $l_{22}$ are coplanar, and provide an Eckardt point$[1: 0: 0: 0]$.

\item $l_2, l_{20}$ and $l_{26}$ are coplanar, and provide an Eckardt point $[0: 1: 0: 0]$.
\item $l_4, l_{23}$ and $l_{27}$ are coplanar, and provide an Eckardt point $[0: \alpha: 1: 1]$.
\item $l_6, l_{19}$ and $l_{26}$ are coplanar, and provide an Eckardt point $[0: \alpha ^2: 1: 1]$.
\item $l_7, l_{18}$ and $l_{27}$ are coplanar, and provide an Eckardt point $[0: \alpha ^5: \alpha ^4: 1]$.
\item $l_8, l_{21}$ and $l_{26}$ are coplanar, and provide an Eckardt point $[0: 1 : \alpha ^5: 1]$.
\item $l_9, l_{12}$ and $l_{26}$ are coplanar, and provide an Eckardt point $[0: \alpha ^3: \alpha: 1]$.
\item $l_{15}, l_{22}$ and $l_{25}$ are coplanar, and provide an Eckardt point $[0: 1 : \alpha ^3: \alpha]$.

\item $l_7, l_9$ and $l_{15}$ are coplanar to each other, but $l_7$ intersects $l_9$ at $[\alpha^3: \alpha^4: 1: \alpha^2]$ and $l_7$ intersects $l_{15}$ at $[\alpha^5: \alpha^2: 1: \alpha]$. They don't intersect at the same point, hence do not provide an Eckardt point.
\item There are other triples that are coplanar to each other but intersect at different points, just like the triple \{$l_7, l_9$ and $l_{15}\}$:

   $ \{l_1, l_{24}, l_{27} \}, \{ l_1, l_{15}, l_{19} \}, \{l_1, l_5, l_8\}, \{l_1, l_{12}, l_{17}\}, \{l_2, l_4, l_{10}\},$

  $ \{l_2, l_{12}, l_{11}\}, \{l_3, l_9, l_{23}\}, \{l_3, l_{14}, l_{18}\},\{l_3, l_{18}, l_{21}\}, \{l_3, l_{6}, l_{16}\},$

$ \{l_5, l_9, l_{11}\}, \{l_5, l_{11}, l_{13}\}, \{l_5, l_{18}, l_{20}\}, \{l_6, l_7, l_{17}\}, \{l_6, l_{24}, l_{22}\},$

$\{l_8, l_{22}, l_{23}\},  \{l_9, l_{10}, l_{24}\}, \{l_{10}, l_{18}, l_{19}\}, \{l_{10}, l_{17}, l_{25}\}, \{l_{11}, l_{16}, l_{27}\},$

$ \{l_{11}, l_{17}, l_{21}\}, \{l_{12}, l_{14}, l_{18}\}, \{l_{12}, l_{22}, l_{24}\}, \{l_{15}, l_{16}, l_{20}\}, \{l_{17}, l_{20}, l_{23}\}. $
\end{itemize}

These give us all the 9 Eckardt points, and the other intersection information can be easily derived from a simple computation, so we omit it here.
\end{exam}

 There are other cubic surfaces that we can compute all the Eckardt points on them.

\begin{exam}
 Clebsch gave a model of a special kind of cubic surfaces, called the Clebsch diagonal surface, where all 27 lines are defined over the field $\mathbb Q[\sqrt{5}]$, and in particular are all real. In addition, up to isomorphism, the Clebsch surface is the only cubic surface that has ten Eckardt points \cite{Hunt}.

The Clebsch diagonal surface in $\bP^3[x,y,z,w]$ is defined by $$ 0= x^3 + y^3 + z^3 + w^3 - (x+y+z+w)^3. $$ We can see the ten Eckardt points on it from the following graph \\
(cited from http://enriques.mathematik.uni-mainz.de/csh/playing/galery/index.html).
\begin{figure}[h!]
  \caption{The Clebsch diagonal surface}
  \centering
    \includegraphics[width=0.3\textwidth]{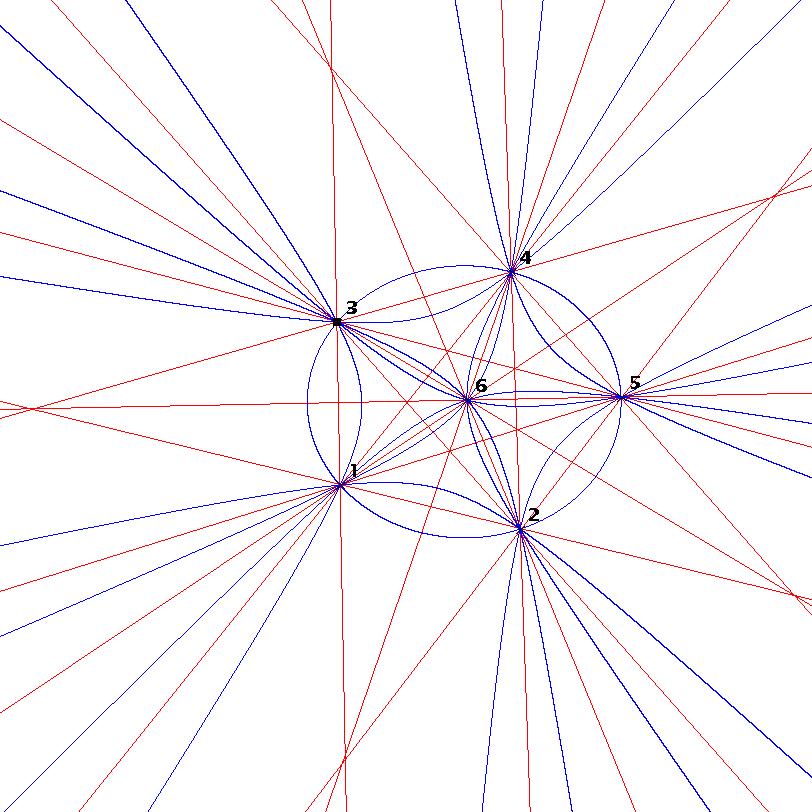}
\end{figure}

\end{exam}
\begin{exam}
The Fermat cubic $$ 0= x^3 + y^3 + z^3 + w^3$$ is the unique one with
18 Eckardt points \cite{Hunt}.
All other cubic surfaces have 0, 1, 2, 3, 4, 6 or 9 such points, as
Beniamino Segre proved in 1942 \cite{Segre}.
\end{exam}
\vspace{1cm}

There exists a cubic surface over a number field $k$ such that the
image of G=$Gal( \bar k /k)$ coincides with the whole Weyl group
$W(E_6)$, i.e.,  the order of the group is 51840. When the field $k$
is $\mathbb Q$, an example is given by Shioda, Tetsuji as follows
\cite{Shioda}.

\begin{exam} The minimum field extension of $\mathbb
Q$ where all the 27 lines of the cubic surface
$$y^2 + 2yz = x^3 + x + xz^2 + z + z^2 + 1$$
are defined is the splitting field of a polynomial of degree 27. The
order of the Galois group is 51840.
\end{exam}

\section{Summary of unirationality results}

However, when our field $k$ is a finite field of small cardinality,
there are some very important geometric properties of the cubic
surface can be derived from the configuration of the 27 lines
defined over the field $\bar k$.

 The following theorem provides
unirationality of all del Pezzo surfaces of degree bigger or equal than 2(Theorem 29.4 \cite{Manin}).

\begin{thm}\label{unirational}
Let X be a del Pezzo surface surface of degree $9-r$ over a number field
$k$ and for $r \geqslant 5$ suppose that there is a point in $X(k)$
which is not on any exceptional curve, and additionally, not on the branch locus when the degree is 2. Then there exists a rational
map $\varphi : \bP^2 \dashrightarrow X$ of degree $\delta _r$ which
is given by the following table:
\begin{center}
\begin{tabular}{lllll}
\hline
r &\vline  $ \ \ r\leqslant 4$ & 5 & 6 & 7 \\
\hline
$\delta _r$ &\vline \ \  1 & 2 & 6 & 24 \\
\hline
\end{tabular}
\end{center}
\end{thm}
Especially, for $r=5$ or 6, i.e., degree 4 del Pezzo surfaces and
cubic surfaces, the requirement that there exists a $k$-point
outside the exceptional curves can be omitted. Therefore, when we
consider cubic surface X over a finite field $k$, the entire
requirement can be omitted, since all finite fields are $C_1$
fields, which is also observed by Koll\'ar's theorem in chapter one.

\subsection{A unirational but birationally non-trivial cubic surface}

  Manin asked a question in Remark \cite{Manin} [Remark 30.1.1] that whether the following cubic surface is unirational.
As an application of the previous theorem, let's answer it in the following example.

\begin{exam}

 Let $k$ be the field of four elements, $\theta \in k$, $\theta \neq 0, 1.$

 Our cubic surface X in $\mathbb P(X,Y,Z,W)$ is given by the equation
 $$X^3+ Y^3+ Z^3 + \theta W^3=0.$$

 For all elements $t \in F_4$, we have $t^3=0$ or 1. Obviously $\theta$ cannot be written as a sum of three cubes. Therefore, all
 the $\bF_4$-points on X lie on the plane section $W=0$. This is the
 non-singular cubic curve $X^3+ Y^3+ Z^3=0$. All $\bF_4$-points on X
 can be easily computed as follows: $[1,1,0,0]$, $[1,\theta,0,0]$,
 $[1,\theta+1,0,0]$, [0,1,1,0], $[0,1,\theta,0], [0,1,\theta+1,0],
 [1,0,1,0], [1,0,\theta,0], [1,0,\theta+1,0]$.

Moreover, X is minimal and birationally non-trivial, since
$$ H^1(Gal(\bar k/k), Pic(X \otimes \bar k))\cong \mathbb Z_3 \times
\mathbb Z_3.$$

Applying Kollar's theorems here, we have that X is unirational over $k$,
and the degree of the unirational map is divisible by 3.

\end{exam}

\section{The rationality of cubic surfaces}

As we have stated before, the rationality of a cubic surface is
determined by the configuration of 27 lines implicitly as long as there is a rational point. We will
study the detailed geometric properties of the cubic surface by
looking at the Galois action on its automorphism group and other related
invariants.

Let we start by setting up the notations. Throughout this section,
we will use X to denote our cubic surface, and we will consider the
surface X over a finite field $k$.

 Following Segre, we denote by $S_n$ any subset of the 27 lines on X that satisfies the conditions below:
 \begin{enumerate}
  \item $S_n$ consists of n lines, no two of which meet.
  \item If $S_n$ contains a line L, then $S_n$ also contains all the conjugates of L over $k$.
  \end{enumerate}
  Because of 1, we have that $n \leq 6$.
\begin{thm}\label{SD}
Let X be a smooth cubic surface over a field $k$.
A necessary and sufficient condition that X should be rational over $k$ is that X should contain a
point defined over $k$ and that X should have at least one $S_2$,
$S_3$ or $S_6$.
\end{thm}

Therefore, we can determine the rationality of cubic surface by
looking at each conjugacy class in $W(E_6)$, and calculate whether
it contains a $S_2$, $S_3$ or $S_6$.

  As we stated before, the Frobenius endomorphism F induces a
  permutation $F^*$ on the 27 lines. The group of those permutations of 27 lines preserves
  incidence relations, and the rationality properties of the lines
  on X depend only on the conjugacy classes in the Weyl group that act on G to which $F^*$
  belongs. These classes have been enumerated by Frame
  (\cite{Frame}; see also table 1 in \cite{Manin}).

\begin{defn}
Suppose that $F$ is an element of the group $Gal(\bar k/k)$ acting
on the Picard group of X. Let $K \supset k$ be the field of
invariants of the cyclic subgroup $(F^n)$. Then the \textbf{index of F} is
the maximal number of geometric components of a divisor on $X\otimes
K$ which can be collapsed over $k$. In particular, if the index of
$F$, $i(F)= 0$, then the surface $X\otimes K$ is minimal.
\end{defn}
\begin{rem}
 The index is given by considering the action of the Weyl group
 $W(E_6)$ on the set of exceptional lines. If we realize the
 corresponding element as an automorphism of a cyclic field
 extension $K/k$ that acts on the Picard group, then the index is the
 maximal number of lines on the cubic surface that can be collapsed over
 $k$.
\end{rem}
\begin{rem}
More concretely, when we consider our surface over the finite field $\bF_q$, the index is the maximal ``n" in the definition of $S_n$. Hence, there is a close relation between the index and the rationality of a cubic surface over the finite fields.
\end{rem}

For example, the theorem given by Manin (Theorem 33.1 in \cite{Manin})
provides some relation between the index of X and the rationality of
X.
\begin{thm}
Every minimal cubic surface over $k$ is birationally non-trivial,
when $k$ is a perfect field.
\end{thm}

  The elements of the Weyl group $W(E_6)$ fall into twenty-five conjugacy classes, which
  Frame calls $C_1$ to $C_{25}$. The following table contains 5 items of information for each
  class:
  \begin{enumerate}
  \item the index of each element in the conjugacy classes.
  \item the order of the element.
  \item trace on Picard group.
  \item the first cohomology group $H^1(G,Pic(X \otimes \bar k))$, where $G$ is the
  cyclic group in $W(E_6)$ generated by the element.
  \item its type as a permutation of the 27 lines.
\end{enumerate}
  These items are explicitly listed in \cite{Manin}. However, there are
  two mistakes on the information of $H^1$ of classes $C_{20}$ and $C_{4}$, which were corrected by
  Urabe in Corollary 1.17, specified with $r=6$\cite{Urabe}. He proved that the group $H^1$ must be a square. We changed the $H^1$ of
classes $C_{20}$ and $C_{4}$ to be 0 as shown in the following table.
\newpage
\begin{table}
\caption{\sl{conjugacy classes in G}}
\begin{tabular}{l l l l l l}
 \hline
  Class & Index   & Order & Trace on Pic & $H^1$ & Permutation type  \\
 \hline
  $C_{13}$ & 0 & 12 & 0 & 0 & ${3, {12}^2}$   \\ [-2ex]
  $C_{12}$ & 0  & 6 & 2 & 0 & ${3, 6^4}$ \\ [-2ex]
  $C_{11}$ & 0  &3  & -2 & $\mathbb Z_3 \times \mathbb Z_3$ & ${3^9}$   \\ [-2ex]
  $C_{14}$& 0 & 9& 1 & 0 & ${9^3}$ \\ [-2ex]
  $C_{10}$ & 0 & 6 & -1 & $\mathbb Z_2 \times \mathbb Z_2$ & ${3, 6^3, 6}$ \\[-2ex]
 $C_{24}$& 1 & 12 & 2 & 0 & {1, 4, 4, 6, 12} \\[-2ex]
 $C_{20}$& 1 & 8  & 1 & 0 & ${1, 2, 8^2, 8}$ \\[-2ex]
 $C_{7}$& 1 & 6   & 2 & 0 & ${1^3, 2^3, 6^3}$ \\[-2ex]
 $C_{19}$ & 1 & 4  & -1 & $\mathbb Z_2 \times \mathbb Z_2$ &${1, 2^3, 4^4, 4}$ \\[-2ex]
 $C_{4}$& 1 & 4   & 3 & 0 & ${1^3,4^6}$ \\[-2ex]
 $C_{3}$& 1 & 2    & -1 & $\mathbb Z_2 \times \mathbb Z_2 $& ${1^3, 2^{12}}$ \\[-2ex]
$C_{25}$& 2 & 10 & 0 & 0 & ${2, 5^2, 5,10}$ \\[-2ex]
 $C_{22}$ & 3 & 6  & -1 &0 & ${3^3, 3^2, 6^2}$ \\[-2ex]
 $C_{8}$& 3 & 6   & 0 & 0 & ${1, 2^2, 2^2, 3^2, 6^2}$ \\[-2ex]
 $C_{23}$& 6 & 6  & 1 & 0 & ${3, 6^2, 6, 6}$ \\[-2ex]
 $C_{15}$& 6 & 5  & 2 & 0 & ${1^2, 5^2, 5^2, 5}$ \\[-2ex]
 $C_{5}$& 6 & 4   & 1 & 0 & ${1, 2^2, 2, 4^2, 4^2, 4}$ \\[-2ex]
 $C_{9}$ & 6 & 3   & 1 & 0 & ${3^6, 3^3}$ \\[-2ex]
 $C_{18}$& 6 & 4  & 3 & 0 & ${1^5, 2, 4^4, 4}$ \\[-2ex]
 $C_{21}$ & 6 & 6  & 2 & 0 & ${1^3, 2^3, 3^4, 6}$\\[-2ex]
 $C_{17}$ & 6 & 2  & 1 & 0 & ${1^3, 2^6, 2^6}$ \\[-2ex]
 $C_{6}$& 6 & 3   & 4 & 0 & ${1^9, 3^6}$ \\[-2ex]
 $C_{2}$& 6 & 2   & 3 & 0 & ${1^7, 2^8, 2}$ \\[-2ex]
 $C_{16}$& 6 & 2  & 5 & 0 & ${1^15, 2^6}$ \\[-2ex]
 $C_{1}$& 6 & 1   & 7 & 0 & ${1^{27}}$ \\
 \hline
\end{tabular}
\end{table}
\newpage

\begin{rem}
The group $H^1$ is a birational invariant. Hence, $H^1 =0$ when the
surface X is birationally trivial, but the converse is not always
true.

\end{rem}

Then what's the relation between the value of the first cohomology
group and the surface with the corresponding associated conjugacy
class? Let's start by analyzing  the 27 lines.

 As in Frame \cite{Frame} and \cite{Sw1}, we denote the 27 lines by those 27 ordered triplets
  made up from the symbols 0, 1, 2, 3 which contain just one 0. Here the
  symbols 0, 1, 2, 3 are to be regarded as elements of the non-cyclic
  group of order four, so that $$1+1=2+2=3+3=0=1+2+3$$ and so on.

  In this notation, two lines are coplanar if and only if the sum of
  these representing triplets is a triplet which represents a
  line (that is to say, which contains just one 0); and if so, this is
  the remaining line in the plane. To make it easier to follow the
  geometrical arguments below, I list the representatives of
  each conjugacy class here (except for class $C_1$, which contains only
  the identity)\cite{Sw2}. These representatives are written as permutations of
  the 27 lines.
\newpage

\begin{center}
\begin{longtable}{ l }
\caption[Representative elements of conjugacy
classes]{Representative elements of conjugacy classes}\\ [-3ex]
 \hline \hline
\endfirsthead
{\tablename} \thetable{} -- Continued \\
\hline
\endhead
\hline
{Continued on Next Page\ldots}\\
\endfoot
\\[-1.8ex] \hline \hline
\endlastfoot
 $C_2$: (011, 022) (013, 301) (023, 302) (031, 310) (032, 320) (101, 220) (102,210) \\[-2ex]
               (201, 120) (202, 110) (303, 330)\\[-2ex]
 $C_3$: (011,022) (013, 032) (023, 031) (101, 120) (102, 110) (103, 130) (201, 220) \\[-2ex]
        (202, 210) (203, 230) (301, 320) (302,310) (303, 330)\\[-2ex]
 $C_4$: (011, 320, 101, 230) (012, 013, 102, 103) (021, 203, 201, 023) (022, 130, \\[-2ex]
        202, 310) (031, 032, 301, 302) (033, 210, 303, 120)\\[-2ex]
 $C_5$: (011, 022) (103, 230) (203, 130) (012, 330, 021, 303) (013, 310, 302, 032)\\[-2ex]
        (023, 031, 301, 320) (101, 102, 210, 220) (201, 202, 110, 120)\\[-2ex]
 $C_6$: (011, 013, 012) (021, 023, 022) (031, 033, 032) (101, 103, 102) (201, 203, 202) \\[-2ex]
        (301, 303, 302)\\[-2ex]
 $C_7$: (011, 101) (022, 202) (033, 303) (012, 013, 230, 102, 320) (021, 310, 203,\\[-2ex]
        201, 130, 023) (031, 032, 120, 301, 302, 210)\\[-2ex]
 $C_8$: (011, 022) (102, 210) (201, 120 ) (303, 330) (012, 130, 203) (021, 103, 230)\\[-2ex]
        (013, 110, 320, 301, 202, 032) (023, 031, 101, 302, 310, 220)\\[-2ex]
 $C_9$: (011, 033, 022) (012, 031, 023) (013, 032, 021) (101, 103, 102) (201, 203, 202) \\[-2ex]
        (301, 303, 302) (110, 130, 120) (210, 230, 220) (310, 330, 320)\\[-2ex]
 $C_{10}$: (101, 202, 303) (011, 022, 102, 110, 220, 201) (012, 330, 021, 210, 033, 120)\\[-2ex]
        (013, 203, 230, 310, 302, 032) (023, 031, 301, 320, 130, 103)\\[-2ex]
 $C_{11}$: (011, 033, 022) (012, 031, 023) (013, 032, 021) (101, 303, 202) (102, 301, 203)\\[-2ex]
        (103, 302, 201) (110, 330, 220) (120, 310, 230) (130, 320, 210)\\[-2ex]
 $C_{12}$: (101, 202, 303) (011, 220, 033, 110, 022, 330) (012, 013, 302, 210, 310, 203)\\[-2ex]
        (021, 103, 320, 120, 301, 023) (031, 032, 201, 130, 230, 102)\\[-2ex]
 $C_{13}$: (110, 220, 330) (011, 022, 302, 103, 023, 031, 101, 202, 032, 013, 203, 301)\\[-2ex]
        (012, 130, 303, 230, 021, 210 ,102, 310 ,033, 320, 201, 120)\\[-2ex]
 $C_{14}$: (011, 022, 202, 032, 013, 303, 023, 031, 101) (012, 130, 110, 033, 220, 230, \\[-2ex]
        021, 310, 320) (102, 103, 330, 203, 201, 120, 301, 302, 210)\\[-2ex]
 $C_{15}$: (011, 032, 330, 031, 022) (012, 320, 310, 021, 303) (013, 033, 023, 302, 301)\\[-2ex]
        (101, 102, 210, 230, 220) (201, 202, 110, 130, 120)\\[-2ex]
 $C_{16}$: (101, 201)(102, 202)(103, 203)(110, 210)(120, 220)(130, 230)\\[-2ex]
 $C_{17}$: (011, 022)(012, 021)(013, 023)(031, 032)(101, 202)(102, 201)(103, 203)(301, 302)\\[-2ex]
        (110, 220)(120, 210)(130, 230)(310, 320)\\[-2ex]
 $C_{18}$: (011, 022)(013, 101, 320, 201)(023, 102, 310, 202)(031, 110, 302, 210)\\[-2ex]
        (032, 120, 301, 220)(103, 203, 130, 230)\\[-2ex]
 $C_{19}$: (011, 022)(103, 130)(203, 230)(012, 330, 021, 303)(013, 310, 302, 032)\\[-2ex]
        (023, 031, 301, 320)(101, 202, 210, 120)(102, 110, 220, 201)\\[-2ex]
 $C_{20}$: (011, 022)(012, 130, 303, 230, 021, 103, 330, 203)(013, 110, 220, 023, 031,\\[-2ex] 101, 202, 032)
        (102, 310, 320, 201, 120, 301, 302, 210)\\[-2ex]
 $C_{21}$: (101, 202)(201, 202)(301, 302)(013, 033, 023)(110, 130, 120)(201, 230, 220)\\[-2ex](310, 330, 320)
        (011, 032, 021, 012, 031, 022)\\[-2ex]
 $C_{22}$: (013, 203, 210)(023, 120, 103)(033, 320, 310)(303, 230, 130)(110, 220, 330)\\[-2ex]
        (011, 022, 302, 012, 021, 301)(031, 101, 202, 032, 102, 201)\\[-2ex]
 $C_{23}$: (011, 022, 033)(012, 110, 220, 021, 101, 202)(013, 130,
 301, 302, 230, 012)\\[-2ex] (031, 103, 310, 320, 203, 032)(012, 330, 201,
 120, 303, 210)\\[-2ex]
  $C_{24}$: (011, 101, 202, 022)(033, 110, 303, 220)(031, 032, 120,
  301, 302, 210)\\[-2ex](012, 013, 023, 021, 103, 320, 201, 130, 230, 102,
  310, 203)\\[-2ex]
  $C_{25}$: (301, 302)(013, 303, 023, 201, 110)(033, 120, 130, 230,
  210)(103, 310, 320, 203, 330)\\[-2ex](011, 022, 102, 032, 202, 012, 021,
  101, 031, 201)\\[-2ex]
\end{longtable}
\end{center}

\begin{thm}
Over a finite field $k=\bF_q$, the rationality of a cubic surface is
governed by the vanishing of all the first cohomology groups
$H^1(\left< \sigma ^n\right>, Pic(X \otimes \bar k))$, where $\sigma$ is the conjugacy class in the
Weyl group $W(E_6)$ associated with the Frobenius and n varies over
all powers.
\end{thm}
\begin{proof}

  Let's first consider the situation that the cubic surfaces are minimal: by
  definition, their associated conjugacy classes are of index 0. Hence there are five classes $C_{13}$, $C_{12}$,
$C_{11}$, $C_{14}$ and $C_{10}$ by Table 1.

The theorem above indicates that they must be birationally
non-trivial over $\bar k$, therefore we need to show that \textbf{to
each of these 5 classes, there exist a power $n \geq 0$ to the class
such that the corresponding first cohomology group $H^1(\left< \sigma ^n\right>,
Pic(X \otimes \bar k))$ is non-vanishing}, where $\sigma=C_{13}$, $C_{12}$, $C_{11}$,
$C_{14}$ and $C_{10}$.

We abuse the notation $\left< \sigma ^n\right>$ here, which I use to represent
the power of the action $G=Gal(\bar k/k)$ on $Pic(X \otimes \bar k)$. Since the
representation of $G$ in $Pic(X \otimes \bar k)$ is given by a finite field extension
of $k$, the image of $G$ in Aut $Pic(X \otimes \bar k)$, which corresponds to a
conjugacy class in $C_1$ to $C_{25}$. Hence the power of the action
is actually the power of the conjugacy class, or in other words, the
power of a representation of each conjugacy class. Therefore, we
only need to show that to each of the 5 classes, there exist a power
of its representation, such that \textbf{the corresponding first
cohomology group $H^1$ is not zero}.

\begin{enumerate}
\item For the class $C_{13}$, take the 4th power of its
representation. The orbit decomposition $\{3, {12}^2\}$ becomes
$\{3, {3}^8\}$ or $\{3^9 \}$. And the representation becomes
\begin{align}\nonumber
 C_{13}^4: &(110, 220, 330)(011, 023, 032)(022, 031, 013)(302, 101,
203)(103, 202, 301)\\ \nonumber &(012, 021, 033)(130, 210,
320,)(303, 102, 201)(230, 310, 120).
\end{align}

Since $110+220=330$, $220+330=110$, $110+330=220$, the first orbit
in the new representation contains three lines which are pairwise
coplanar. Similarly, we have 9 orbits of three disjoint lines,
hence the decomposition type of $ C_{13}^4$ is $\{3^9\}$. Checking table 2, we can see it is the same as $C_{11}$. We can
easily see from table 1 that the value of its first cohomology group
is $\mathbb Z_3 \times \mathbb Z_3$, not zero.

\item For the class $C_{11}$, the value of its first cohomology group is $\mathbb
Z_3 \times \mathbb Z_3$, not zero.

\item For the class $C_{12}$, take the 2nd power of its
representation, which becomes
\begin{align}\nonumber
 C_{12}^2:&(101, 202, 303)(011, 033, 022)(220, 110, 330)(012,
302, 310) \\\nonumber& (013, 210, 203)(021, 320, 301)(103, 120,
023)(031, 201, 230)(032, 130, 120).
\end{align}

It is the same as $C_{11}$ again. So the value of its first
cohomology group is $\mathbb Z_3 \times \mathbb Z_3$, not zero.

\item For the class $C_{14}$, take the 3rd power of its
representation, which becomes
\begin{align}\nonumber
C_{14}^3: &(011, 032, 023)(022, 013, 031)(202, 303, 101)(012, 033,
021)\\\nonumber &(130, 220, 310)(110, 230, 320)(102, 203, 301)(103,
201, 302)(330, 120, 210).
\end{align}

It is the same as class $C_{11}$, hence the corresponding first
cohomology is not zero.
\item For the class $C_{10}$, the value of its first cohomology group is $\mathbb Z_2 \times \mathbb Z_2$, not zero.
\end{enumerate}

\vspace{0.7cm}
  Second, let's consider the classes that contain a pair of
disjoint lines. By the Theorem \ref{SD}, if the surface X is represented
by one of these classes, X is birationally trivial, hence the value of the first cohomology groups of every power are zero.

There are 7 of them by table 1:
\begin{itemize}
\item $C_{5}: (103, 230)$ is a pair of disjoint line since $1+2=3,
0+3=3, 3+0=3$, i.e. there is no single zero in the sum.
\item $C_{16}: (102, 202)$, $C_{17}: (103, 203)$, $C_{8}: (102, 210)$, $C_{21}: (201, 202)$, $C_{2}: (013, 301)$ contain a pair of disjoint lines, similarly as
$C_{5}$.
\end{itemize}

\vspace{0.8
cm}
  Let's consider the classes that contain three disjoint lines. By Theorem \ref{SD}, if the surface X is represented
by one of these classes, X is
birationally trivial, hence the value of the first cohomology groups of every power are zero. There are 3 of them by the table:
\begin{itemize}
\item $C_{6}: (011, 013, 012)$ in the representation is an orbit of 3 disjoint lines since $011+013=002$,
$011+012=003$, $012+013=001$, i.e. these three lines are pairwise
disjoint.
\item $C_{9}: (101, 103, 102)$ and $C_{22}: (033, 320, 310)$ in the representation is an orbit of 3 disjoint lines.
\end{itemize}

 \vspace{1cm}
 Let's consider the classes that contain four disjoint lines:
 \begin{itemize}
   \item $C_{18}:(013, 101, 320, 201)$
 \end{itemize}

     The four disjoint lines can be blown down, where the surface X becomes $Bl_{\{p_1,p_2\}} \bP^2$, in
other words, a surface contains two disjoint lines, therefore by
Theorem \ref{SD}, it is rational, so is X.

Therefore, the value of its first cohomology groups of every power are zero.

\vspace{1cm}

 Let's consider the class that contains five disjoint lines, which is $C_{15}:(013, 033, 023, 302, 301)$. Since the index of $C_{15}$ is 6, by the Theorem \ref{SD}, it is rational.

 Hence, all $H^1$ groups are zero.

\vspace{1cm}
 Let's consider the conjugacy class that contain an orbit of six disjoint lines, which is $C_{23}:(013, 130, 301, 302, 230, 023)$.
 Similar as before, the value of all the first cohomology groups are zero.

\vspace{1cm}

  Let's consider all the remaining classes one by one:

\begin{itemize}

\item For $C_{7}$, it contains 3 orbits of 2 coplanar lines, 3 orbits of 6 lines. These three orbits
are as follows,
\begin{enumerate}
\item (012, 013, 230, 102, 103, 320)
\item (021, 310, 203, 201, 130, 023)
\item (031, 032, 120, 301, 302, 210).
\end{enumerate}
In these three orbits, each line is coplanar with another line, and
disjoint from all the other 5 lines. It doesn't contain any $S_2$,
$S_3$ or $S_6$, hence by Theorem \ref{SD} it is not birationally
trivial. Let's look at the value of its first cohomology groups.

 We have the relation $C_{7}^2=C_{6}$,\emph{$C_{7}^3=C_{3}$}, $C_{7}^6=$identity, and all the other higher
powers of $C_{7}$ are itself. By the information in table 1, we can
conclude that the values of all the $H^1$ groups are zero except for the
third power of its representation.

 $C_{7}$ is not rational although its first cohomology group vanishes.

\item For $C_{24}$, it has an orbit of 4 coplanar lines, an orbit of 2
pairs of 2 coplanar lines, an orbit of 3 pairs of 2 coplanar
lines, and an orbit of 12 lines. None of them are skew or
coplanar. By Theorem \ref{SD} and the analysis in the previous steps, it
is not birationally trivial. Let's look at the values of its first
cohomology groups.

We have the relations $C_{24}^2=C_{8}$, $C_{24}^3=C_{19}$,
$C_{24}^4=C_{21}$, and all power bigger than 6 except 12 of $C_{24}$
are itself. By the information in table 1, we can conclude that the
value of all the first cohomology groups are not all vanishing
because $H^1(C_{24}^3)=H^1(C_{19}) \neq 0$.
\item For $C_{20}$, it has one orbit of 2 coplanar lines, and 3 orbits of
8 lines. Hence by the Theorem \ref{SD}, it is not rational. By
the information in table 1, the value of its first cohomology group is zero. However,
$C_{20}^2=C_{4}$, $C_{20}^4=C_{3}$, so $H^1(C_{20}^4)=H^1(C_{4}) \neq 0$.

 $C_{20}$ is not rational although its first cohomology group vanishes.

\item For $C_{19}$, it has three orbits of 2 coplanar lines, and 5 orbits
of 4 coplanar lines. Hence by the Theorem \ref{SD}, it is not
rational. By the information in, the value of its first
cohomology group is not zero.

\item For $C_{4}$, it has 6 orbits of 4 coplanar lines. Hence by the
Theorem \ref{SD}, it is not rational. By information of table 1, the
value of its first cohomology group is not zero.

\item For $C_{3}$, it has 12 orbits of 2 coplanar lines. Hence by Theorem \ref{SD}, it is not rational. By the information in table 1, the
value of its first cohomology group is not zero.

\end{itemize}
\end{proof}

\section{Rationality of other del Pezzo surfaces}
\subsection{Rationality of del Pezzo surfaces of higher degrees}

For del Pezzo surface of higher degrees, their rationality criteria are given more explicitly. The reason why these are easier than dP1,
dP2 and cubic surface is they have much clearer geometric structures.

\begin{thm}
Let X be a del Pezzo surface of degree $9-r$ and let there exist a
rational map of finite degree $\varphi : \bP^2 \dashrightarrow X$.
Then the degree of $\varphi$ is divisible by the least common
multiple $d_r$ of the exponents of the groups $H^1(G,N_r)$ for all
possible subgroups $G \subset W(R_r)$, where $N_r$ is the lattice
generated by the vectors $l_0, l_1, \ldots, l_r$ in the standard
realization of $R_r$.

The number $d_r$ are given by the following table:
\begin{center}
\begin{tabular}{llllll}
\hline
r &\vline  $r\leqslant 4$ & 5 & 6 & 7 &8\\
\hline
$\delta _r$ &\vline  1 & 2 & 6 & 24 & 120\\
\hline
\end{tabular}
\end{center}
\end{thm}
\begin{proof}
This is Theorem 29.2 and Theorem 29.3 in \cite{Manin}, and the above
table is a generalization of the table in Theorem \ref{unirational}.
\end{proof}

 We can easily see that for del Pezzo surfaces of
degree greater than or equal to 4, their rationality is equivalent to
the existence of a rational point.

\subsection{Rationality of del Pezzo surfaces of degree 2}

 In the previous section, we have shown that the rationality of a cubic
 surface over a finite field is entirely determined by the configuration of the 27
 lines on it, or more precisely the vanishing of all the first cohomology groups.

 Comparing the cubic surfaces with the del Pezzo surfaces of other degrees,
 we know that they are generated in a very similar way, which is by blowing up $\bP^2$ along different number of points in general
 position. Also we know that all the lines (i.e. exceptional curves) on the del Pezzo surfaces are
 generated in a very similar way. Therefore, we naturally expect that the
 rationality criterion is similar for cubic surfaces and other del
 Pezzo surfaces.

 But it is NOT always true. Let's look at the following examples:

\begin{exam}
There exists a minimal del Pezzo surface of degree 2 over a finite field which is
birationally non-trivial but has all vanishing first cohomology groups.

Let S denotes the del Pezzo surface of degree 2 which is represented
by Carter symbol $(4A_1)'$ (\cite{Urabe} or \cite{Sw2}), with index
0, hence minimal.

By the table of classes of conjugacy elements in the Weyl group
$W(E_7)$ given by Urabe, we can see that the value of its first
cohomology group is 0, and the orbit decomposition is $\{2^4, 2^{24}
\}$. Therefore, taking the second power it becomes the identity, with the first cohomology group vanishing.

In conclusion, S is birationally non-trivial with all the first cohomology groups to any power vanishing.
\end{exam}

Checking other minimal del Pezzo surfaces of degree 2, we can see that S is not the only minimal dP2 with all vanishing first cohomology groups. We list all such dP2s as follows:
\begin{itemize}
\item The conjugacy class $C_5$ with Carter symbol $D_6$ and orbit decomposition \\ $\{2,2^2,10,10,10,{10}^2 \}$.
\item The conjugacy class $C_{16}$ with Carter symbol $E_7(a_1)$ and orbit decomposition $\{14, 14, 14, 14 \}$.
\end{itemize}
All the other  16 classes of conjugacy elements in the Weyl group corresponding to the minimal dP2s have at least one of its first cohomology groups non-trivial.

\subsection{Rationality of del Pezzo surfaces of degree 1}
We have two different situations with the cubic surfaces and the dP2s when we consider the rationality of them. How about the rationality of del Pezzo surfaces of degree 1? Is it like the cubic surfaces or the dP2s?

We have the following theorem:
\begin{thm}
All minimal del Pezzo surfaces of degree one have a nontrivial first cohomology group to some power of their corresponding conjugacy class.
\end{thm}

\begin{proof}

Let's look at all the conjugacy classes corresponding to the minimal del Pezzo surfaces of degree 1 (dP1). There are 37 classes of conjugacy elements in the Weyl group of $W(E_8)$ according to table 2\cite{Urabe}. And there are 12 classes of them with their first cohomology groups vanishing.

 Let's check them one by one as follows:
\begin{enumerate}
\item For class $C_5$ with Carter symbol $D_7$ and the orbit decomposition \\ $ \{2^2,4^2,12,12,12,12,12,{12}^2 ,{12}^2,{12}^2,{12}^2,{12}^2,{12}^4\}$, we have the relation \\$C_5^3=C_3$, which has nontrivial first cohomology group.
\item For class $C_6$ with Carter symbol $D_7(a_1)$ and the orbit decomposition\\ $ \{4,4^2,4^2,10,10,10,10,{10}^2 ,{20}^2,{20}^2,{20}^2,{20}^2\}$, we have the relation $C_6^4=C_{102}$, which has nontrivial first cohomology group.
\item For class $C_7$ with Carter symbol $D_7(a_2)$ and the orbit decomposition\\ $ \{6,6,6^2,8,8,8,8^2,8^2,8^2, {24}^2 ,{24}^2 ,{24}^2\}$, we have the relation $C_7^3=C_4$, which has nontrivial first cohomology group.
\item For class $C_{29}$ with Carter symbol $E_8$ and the orbit decomposition\\ $ \{30,30,30,30,30,30,30,30 \}$, we have the relation $C_{29}^6=C_{11}$, which has nontrivial first cohomology group.
\item For class $C_{30}$ with Carter symbol $E_8(a_1)$ and the orbit decomposition\\ $ \{24,24,24,24,24,24,24,24,{24}^2 \}$, we have the relation $C_{30}^3=C_{23}$, which has nontrivial first cohomology group.
\item For class $C_{31}$ with Carter symbol $E_8(a_2)$ and the orbit decomposition\\$ \{20,20,20,20,20,20,20,20,{20}^2,{20}^2 \}$, we have the relation $C_{31}^4=C_{11}$, which has nontrivial first cohomology group.
\item For class $C_{32}$ with Carter symbol $E_8(a_3)$ and the orbit decomposition\\$ \{{12}^6,{12}^6,{12}^8 \}$. We have the relation $C_{32}^3=C_{17}$, which has nontrivial first cohomology group.
\item For class $C_{33}$ with Carter symbol $E_8(a_4)$ and the orbit decomposition\\$ \{6,18,18,18,18,{18}^3,{18}^3,{18}^3 \}$, we have the relation $C_{33}^2=C_{14}$, which has nontrivial first cohomology group.
\item For class $C_{34}$ with Carter symbol $E_8(a_5)$ and the orbit decomposition\\$ \{ {15}^2,{15}^2,{15}^2,{15}^2,{15}^2,{15}^2,{15}^2,{15}^2  \}$, we have the relation $C_{34}^3=C_{11}$, which has nontrivial first cohomology group.
\item For class $C_{35}$ with Carter symbol $E_8(a_6)$ and the orbit decomposition\\$ \{ {10}^{12},{10}^{12} \}$, we have the relation $C_{35}^2=C_{11}$, which has nontrivial first cohomology group.
\item For class $C_{36}$ with Carter symbol $E_8(a_7)$ and the orbit decomposition\\$ \{{6}^{4},12,12,{12}^{8},{12}^{8} \}$, we have the relation $C_{36}^2=C_{25}$, which has nontrivial first cohomology group.
\item For class $C_{37}$ with Carter symbol $E_8(a_8)$ and the orbit decomposition\\$ \{{6}^{40}\}$, we have the relation $C_{37}^2=C_{9}$, which has nontrivial first cohomology group.
\end{enumerate}

\end{proof}

\newpage


\newpage
\begin{bibdiv}
\begin{biblist}

 \bib{Frame}{article}{
   author={Frame, J. S.},
   title={The classes and representations of the groups of $27$ lines and
   $28$ bitangents},
   journal={Ann. Mat. Pura Appl. (4)},
   volume={32},
   date={1951},
   pages={83--119},
   issn={0003-4622},
   review={\MR{0047038 (13,817i)}},
}

\bib{GHS}{article}{
   author={Graber, Tom},
   author={Harris, Joe},
   author={Mazur, Barry},
   author={Starr, Jason},
   title={Rational connectivity and sections of families over curves},
   language={English, with English and French summaries},
   journal={Ann. Sci. \'Ecole Norm. Sup. (4)},
   volume={38},
   date={2005},
   number={5},
   pages={671--692},
   issn={0012-9593},
   review={\MR{2195256 (2006j:14044)}},
}
\bib{Hunt}{book}{
    author={Hunt, Bruce},
    title ={The Geometry of Some Special Arithmetic Quotients},
    series= {New York: Springer-Verlag},
    pages ={122--128},
    date= {1996},
}

\bib{Isk79}{article}{
   author={Iskovskih, X. A.},
   title={Minimal models of rational surfaces over arbitrary fields},
   language={Russian},
   journal={Izv. Akad. Nauk SSSR Ser. Mat.},
   volume={43},
   date={1979},
   number={1},
   pages={19--43, 237},
   issn={0373-2436},
   review={\MR{525940 (80m:14021)}},
}

\bib{Kollar1}{book}{
   author={Koll{\'a}r, J{\'a}nos},
   title={Rational curves on algebraic varieties},
   series={Ergebnisse der Mathematik und ihrer Grenzgebiete. 3. Folge. A
   Series of Modern Surveys in Mathematics [Results in Mathematics and
   Related Areas. 3rd Series. A Series of Modern Surveys in Mathematics]},
   volume={32},
   publisher={Springer-Verlag},
   place={Berlin},
   date={1996},
   pages={viii+320},
   isbn={3-540-60168-6},
   review={\MR{1440180 (98c:14001)}},
}

\bib{Kollar02}{article}{
    author={Koll{\'a}r, J{\'a}nos},
   title={Unirationality of Cubic Hypersurfaces},
   journal={Journal of the Institute of Mathematics of Jussieu},
   volume={1},
   pages={467-476 },
   issn={doi:10.1017/S1474748002000117 },

}

\bib{Kollar}{article}{
   author={Koll{\'a}r, J{\'a}nos},
   title={Looking for rational curves on cubic hypersurfaces},
   note={Notes by Ulrich Derenthal},
   conference={
      title={Higher-dimensional geometry over finite fields},
   },
   book={
      series={NATO Sci. Peace Secur. Ser. D Inf. Commun. Secur.},
      volume={16},
      publisher={IOS},
      place={Amsterdam},
   },
   date={2008},
   pages={92--122},
   review={\MR{2484078 (2009k:14047)}},
}
\bib{KMM}{article}{
   author={Koll{\'a}r, J{\'a}nos},
   author={Miyaoka, Yoichi},
   author={Mori, Shigefumi},
   title={Rationally connected varieties},
   journal={J. Algebraic Geom.},
   volume={1},
   date={1992},
   number={3},
   pages={429--448},
   issn={1056-3911},
   review={\MR{1158625 (93i:14014)}},
}

\bib{Manin}{book}{
   author={Manin, Yu. I.},
   title={Cubic forms},
   series={North-Holland Mathematical Library},
   volume={4},
   edition={2},
   note={Algebra, geometry, arithmetic;
   Translated from the Russian by M. Hazewinkel},
   publisher={North-Holland Publishing Co.},
   place={Amsterdam},
   date={1986},
   pages={x+326},
   isbn={0-444-87823-8},
   review={\MR{833513 (87d:11037)}},
}
\bib{Sko01}{book}{
   author={Skorobogatov, Alexei},
   title={Torsors and rational points},
   series={Cambridge Tracts in Mathematics},
   volume={144},
   publisher={Cambridge University Press},
   place={Cambridge},
   date={2001},
   pages={viii+187},
   isbn={0-521-80237-7},
   review={\MR{1845760 (2002d:14032)}},
}

\bib{Segre}{book}{
   author={Segre, B.},
   title={The Non-singular Cubic Surfaces},
   publisher={Oxford University Press},
   place={Oxford},
   date={1942},
   pages={xi+180},
   review={\MR{0008171 (4,254b)}},
}

\bib{Shioda}{article}{
 author= { Tetsuji Shioda},
 title = {The Algebraic Equation of Type E6 and Chebotarev Density},
 conference ={
       title ={ Algebra and Analysis (Proc. Chebotarev centennial conference at Kazan Univ. June 1994 , ed. by Arslanov, Parshin, Shafarevich)},
       publisher ={ Walter de Gruyter & Co.} ,
       },
 pages={109--114},
 date={1996},
}

\bib{Sw1}{article}{
   author={Swinnerton-Dyer, Peter},
   title={Counting points on cubic surfaces. II},
   conference={
      title={Geometric methods in algebra and number theory},
   },
   book={
      series={Progr. Math.},
      volume={235},
      publisher={Birkh\"auser Boston},
      place={Boston, MA},
   },
   date={2005},
   pages={303--309},
   review={\MR{2166089 (2006e:11088)}},

}
\bib{Sw2}{article}{
   author={Swinnerton-Dyer, H. P. F.},
   title={The zeta function of a cubic surface over a finite field},
   journal={Proc. Cambridge Philos. Soc.},
   volume={63},
   date={1967},
   pages={55--71},
   review={\MR{0204414 (34 \#4256)}},
}
\bib{Sw3}{article}{
   author={Swinnerton-Dyer, H. P. F.},
   title={The birationality of cubic surfaces over a given field},
   journal={Michigan Math. J.},
   volume={17},
   date={1970},
   pages={289--295},
   issn={0026-2285},
   review={\MR{0279099 (43 \#4825)}},
}


\bib{Sw4}{article}{
  author={Swinnerton-Dyer, H. P. F.},
   title={Universal equivalence for cubic surfaces over finite and local fields},
   conference={
      title={Symposia Mathematica(Sympos., INDAM, Rome,1979)},
   volume={XXIV },
   },
   book={
       series={Academic Press},
       publisher={London-New York},
       },
   pages={111--143},
   review={\MR{619244 (82k : 14019)}},
}

\bib{Sw5}{article}{
 author={Swinnerton-Dyer, H. P. F.},
   title={Rational points on del Pezzo surfaces of degree 5},
   conference={
   title= {5th Nordic Summer School in Math},
   },
   year ={1970},
   pages ={287--290},
}
\bib{Urabe}{article}{
   author={Urabe, Tohsuke},
   title={Calculation of Manin's invariant for Del Pezzo surfaces},
   journal={Math. Comp.},
   volume={65},
   date={1996},
   number={213},
   pages={247--258, S15--S23},
   issn={0025-5718},
   review={\MR{1322894 (96f:14047)}},
   doi={10.1090/S0025-5718-96-00681-3},
}

\bib{Weil}{article}{
   author={Weil, Andr{\'e}},
   title={Abstract versus classical algebraic geometry},
   conference={
      title={Proceedings of the International Congress of Mathematicians,
      1954, Amsterdam, vol. III},
   },
   book={
      publisher={Erven P. Noordhoff N.X., Groningen},
   },
   date={1956},
   pages={550--558},
   review={\MR{0092196 (19,1078a)}},
}

\bib{Ct}{article}{
 author = {Jean-Louis Colliot-Th\'{e}l\`{e}ne},
 title ={ Points rationnels sur les vari\`{e}t\`{e}s non de type g\`{e}n\`{e}ral},
 journal ={ Cours IHP, lecture notes},
  year ={2003},
}

\bib{Yan}{article}{

 author={ V. I. Yanchevskii},
 title={K-unirationality of conic bundles and splitting fields of simple central algebras},
 journal={ Doklady Akad. Nauk BSSR (Minsk,Belarus) },
 volume= {29},
 date={1985}
 pages ={1061-1064},

}
\bib{Cm}{article}{
  author = {D. F. Coray and M. A. Tsfasman},
 title = {A.: Arithmetic on singular del Pezzo surfaces},
 journal = {Proc. London Math. Soc},
  volume = {3},
   date = {1988},
  pages = {57},
 review={\MR{728739}},
}
\bib{Yu}{article}{
   author={Manin, Yu. I.},
   author={Tsfasman, M. A.},
   title={Rational varieties: algebra, geometry, arithmetic},
   language={Russian},
   journal={Uspekhi Mat. Nauk},
   volume={41},
   date={1986},
   number={2(248)},
   pages={43--94},
   issn={0042-1316},
   review={\MR{842161 (87k:11065)}},
}
\end{biblist}
\end{bibdiv}
\end{document}